\let\noarrow = t
\input eplain

\let\noarrow = t

\input eplain


\magnification=\magstep1

\topskip1truecm
\def\pagewidth#1{
  \hsize=#1
}

\def\pageheight#1{
  \vsize=#1
}

\pageheight{23.5truecm} \pagewidth{16truecm}

\abovedisplayskip=3mm \belowdisplayskip=3mm
\abovedisplayshortskip=0mm \belowdisplayshortskip=2mm
\parindent1pc

\normalbaselineskip=13pt \baselineskip=13pt

\def\spacing{{\smallskip}}

\voffset=0pc \hoffset=0pc


\newdimen\abstractmargin
\abstractmargin=3pc


\newdimen\footnotemargin
\footnotemargin=1pc


\font\eightrm=cmr8 \relax 
\font\sixrm=cmr6 \relax 
\font\eighti=cmmi8 \relax     \skewchar\eighti='177 
\font\sixi=cmmi6 \relax       \skewchar\sixi='177   
\font\eightsy=cmsy8 \relax    \skewchar\eightsy='60 
\font\sixsy=cmsy6 \relax      \skewchar\sixsy='60   
\font\eightbf=cmbx8 \relax 
\font\sixbf=cmbx6 \relax   
\font\eightit=cmti8 \relax 
\font\eightsl=cmsl8 \relax 
\font\eighttt=cmtt8 \relax 

\catcode`\@=11
\newskip\ttglue

\def\eightpoint{\def\rm{\fam0\eightrm}%
 \textfont0=\eightrm \scriptfont0=\sixrm
 \scriptscriptfont0=\fiverm
 \textfont1=\eighti \scriptfont1=\sixi
 \scriptscriptfont0=\fivei
 \textfont2=\eightsy \scriptfont2=\sixsy
 \scriptscriptfont2=\fivesy
 \textfont3=\tenex \scriptfont3=\tenex
 \scriptscriptfont3=\tenex
 \textfont\itfam\eightit \def\it{\fam\itfam\eightit}%
 \textfont\slfam\eightsl \def\sl{\fam\slfam\eightsl}%
 \textfont\ttfam\eighttt \def\tt{\fam\ttfam\eighttt}%
 \textfont\bffam\eightbf \scriptfont\bffam\sixbf
   \scriptscriptfont\bffam\fivebf \def\bf{\fam\bffam\eightit}%
 \tt \ttglue=.5em plus.25em minus.15em
 \normalbaselineskip=9pt
 \setbox\strutbox\hbox{\vrule height7pt depth3pt width0pt}%
 \let\sc=\sixrm \let\big=\eifgtbig \normalbaselines\rm}


 \font\titlefont=cmbx12 scaled\magstep1
 \font\sectionfont=cmbx12
 \font\ssectionfont=cmsl10
 \font\claimfont=cmsl10

 \font\normalfont=cmr10

\catcode`\@=11 \font\teneusm=eusm10 
\font\seveneusm=eusm7  \font\fiveeusm=eusm5
\newfam\eusmfam \textfont\eusmfam=\teneusm
\scriptfont\eusmfam=\seveneusm \scriptscriptfont\eusmfam=\fiveeusm
\def\hexnumber@#1{\ifcase#1
0\or1\or2\or3\or4\or5\or6\or7\or8\or9\or         A\or B\or C\or D\or
E\or F\fi } \edef\eusm@{\hexnumber@\eusmfam}
\def\euscr{\ifmmode\let\next\euscr@\else
\def\next{\errmessage{Use \string\euscr\space only in math mode}}\fi\next}
\def\euscr@#1{{\euscr@@{#1}}} \def\euscr@@#1{\fam\eusmfam#1} \catcode`\@=12

\catcode`\@=11 \font\teneuex=euex10 
 \font\seveneuex=euex7  \newfam\euexfam
\textfont\euexfam=\teneuex  \scriptfont\euexfam=\seveneuex
 \def\hexnumber@#1{\ifcase#1
0\or1\or2\or3\or4\or5\or6\or7\or8\or9\or         A\or B\or C\or D\or
E\or F\fi } \edef\euex@{\hexnumber@\euexfam}
\def\euscrex{\ifmmode\let\next\euscrex@\else
\def\next{\errmessage{Use \string\euscrex\space only in math mode}}\fi\next}
\def\euscrex@#1{{\euscrex@@{#1}}} \def\euscrex@@#1{\fam\euexfam#1}
\catcode`\@=12

\catcode`\@=11 \font\teneufb=eufb10 
\font\seveneufb=eufb7  \font\fiveeufb=eufb5
\newfam\eufbfam \textfont\eufbfam=\teneufb
\scriptfont\eufbfam=\seveneufb \scriptscriptfont\eufbfam=\fiveeufb
\def\hexnumber@#1{\ifcase#1
0\or1\or2\or3\or4\or5\or6\or7\or8\or9\or         A\or B\or C\or D\or
E\or F\fi } \edef\eufb@{\hexnumber@\eufbfam}
\def\euscrfb{\ifmmode\let\next\euscrfb@\else
\def\next{\errmessage{Use \string\euscrfb\space only in math mode}}\fi\next}
\def\euscrfb@#1{{\euscrfb@@{#1}}} \def\euscrfb@@#1{\fam\eufbfam#1}
\catcode`\@=12

\catcode`\@=11 \font\teneufm=eufm10 
\font\seveneufm=eufm7  \font\fiveeufm=eufm5
\newfam\eufmfam \textfont\eufmfam=\teneufm
\scriptfont\eufmfam=\seveneufm \scriptscriptfont\eufmfam=\fiveeufm
\def\hexnumber@#1{\ifcase#1
0\or1\or2\or3\or4\or5\or6\or7\or8\or9\or         A\or B\or C\or D\or
E\or F\fi } \edef\eufm@{\hexnumber@\eufmfam}
\def\euscrfm{\ifmmode\let\next\euscrfm@\else
\def\next{\errmessage{Use \string\euscrfm\space only in math mode}}\fi\next}
\def\euscrfm@#1{{\euscrfm@@{#1}}} \def\euscrfm@@#1{\fam\eufmfam#1}
\catcode`\@=12

\catcode`\@=11 \font\teneusb=eusb10 
\font\seveneusb=eusb7  \font\fiveeusb=eusb5
\newfam\eusbfam \textfont\eusbfam=\teneusb
\scriptfont\eusbfam=\seveneusb \scriptscriptfont\eusbfam=\fiveeusb
\def\hexnumber@#1{\ifcase#1
0\or1\or2\or3\or4\or5\or6\or7\or8\or9\or         A\or B\or C\or D\or
E\or F\fi } \edef\eusb@{\hexnumber@\eusbfam}
\def\euscrsb{\ifmmode\let\next\euscrsb@\else
\def\next{\errmessage{Use \string\euscrsb\space only in math mode}}\fi\next}
\def\euscrsb@#1{{\euscrsb@@{#1}}} \def\euscrsb@@#1{\fam\eusbfam#1}
\catcode`\@=12

\catcode`\@=11 \font\tenmsa=msam10 
\font\sevenmsa=msam7  \font\fivemsa=msam5
\font\tenmsb=msbm10  \font\sevenmsb=msbm7
 \font\fivemsb=msbm5 \newfam\msafam
\newfam\msbfam \textfont\msafam=\tenmsa
\scriptfont\msafam=\sevenmsa
  \scriptscriptfont\msafam=\fivemsa
\textfont\msbfam=\tenmsb  \scriptfont\msbfam=\sevenmsb
  \scriptscriptfont\msbfam=\fivemsb
\def\hexnumber@#1{\ifcase#1 0\or1\or2\or3\or4\or5\or6\or7\or8\or9\or
        A\or B\or C\or D\or E\or F\fi }
\edef\msa@{\hexnumber@\msafam} \edef\msb@{\hexnumber@\msbfam}
\mathchardef\square="0\msa@03 \mathchardef\subsetneq="3\msb@28
\mathchardef\supsetneq="3\msb@29 \mathchardef\ltimes="2\msb@6E
\mathchardef\rtimes="2\msb@6F \mathchardef\dabar="0\msa@39
\mathchardef\daright="0\msa@4B \mathchardef\daleft="0\msa@4C

\def\Bbb{\ifmmode\let\next\Bbb@\else
        \def\next{\errmessage{Use \string\Bbb\space only in math mode}}\fi\next}
\def\Bbb@#1{{\Bbb@@{#1}}}
\def\Bbb@@#1{\fam\msbfam#1}
\catcode`\@=12



\newcount\senu
\def\senum{\number\senu}
\newcount\ssnu
\def\ssnum{\number\ssnu}
\newcount\fonu
\def\fonum{\number\fonu}

\def\num{{\senum.\ssnum}}
\def\numfo{{\senum.\ssnum.\fonum}}


\outer\def\section#1\par{\vskip0pt
  plus.3\vsize\penalty20\vskip0pt
  plus-.3\vsize\bigskip\vskip\parskip
  \message{#1}\centerline{\sectionfont\senum\enspace#1.}
  \nobreak\smallskip}

\def\endsection{\advance\senu by1\penalty-20\smallskip\ssnu=1}
\outer\def\ssection#1\par{\bigskip
  \message{#1}{\noindent\bf\num\ssectionfont\enspace#1.\thinspace}
  \nobreak\normalfont}

\def\endssection{\advance\ssnu by1\smallskip\ifdim\lastskip<\medskipamount
\removelastskip\penalty55\medskip\fi\fonu=1\normalfont}

\def\proclaim #1\par{\bigskip
  \message{#1}{\noindent\bf\num\enspace#1.\thinspace}
  \nobreak\claimfont}

\def\cor{\proclaim Corollary\par}
\def\defi{\proclaim Definition\par}
\def\lemma{\proclaim Lemma\par}
\def\prop{\proclaim Proposition\par}
\def\rmk{\proclaim Remark\par\normalfont}
\def\thm{\proclaim Theorem\par}

\def\endcor{\endssection}
\def\enddefi{\endssection}
\def\endlemma{\endssection}
\def\endprop{\endssection}
\def\endrmk{\endssection}
\def\endthm{\endssection}

\def\Proof{{\noindent\sl Proof: \/}}


\def\maplefto#1{\ \smash{\mathop{\longleftarrow}\limits^{#1}}\ }

\def\llongrightarrow{\relbar\joinrel\relbar\joinrel\rightarrow}
\def\lllongrightarrow{\hbox to 40pt{\rightarrowfill}}

\def\twoheadrightarrow{\rightarrow\kern -8pt\rightarrow}

\def\maprighto#1{\smash{\mathop{\longrightarrow}\limits^{#1}}}

\def\mapdownr#1{\Big\downarrow\rlap{$\vcenter{\hbox{$\scriptstyle#1$}}$}}
\def\mapdownl#1{\llap{$\vcenter{\hbox{$\scriptstyle#1$}}$}\Big\downarrow}

\def\llongmaprighto#1{\ \smash{\mathop{\llongrightarrow}\limits^{#1}}\ }

\def\lllongmaprighto#1{\ \smash{\mathop{\lllongrightarrow}\limits^{#1}}\ }

\def\llongleftarrow{\leftarrow\joinrel\relbar\joinrel\relbar}

\def\longleftmapsto{\longleftarrow\kern-2pt\mapstochar\;}

\def\llongmapsto{\,\vert\kern-3.2pt\joinrel\longrightarrow\,}
\def\llongmapsto{\,\vert\kern-3.7pt\joinrel\llongrightarrow\,}
\def\lllongmapsto{\,\vert\kern-5.5pt\joinrel\lllongrightarrow\,}

\def\isomarrow{\maprighto{\lower3pt\hbox{$\scriptstyle\sim$}}}
\def\llongisomarrow{\llongmaprighto{\lower3pt\hbox{$\scriptstyle\sim$}}}
\def\lllongisomarrow{\lllongmaprighto{\lower3pt\hbox{$\scriptstyle\sim$}}}

\def\lisomarrow{\maplefto{\lower3pt\hbox{$\scriptstyle\sim$}}}

\font\labprffont=cmtt8
\def\strutdepth{\dp\strutbox}
\def\labtekst#1{\vtop to \strutdepth{\baselineskip\strutdepth\vss\llap{{\labprffont #1}}\null}}
\def\marglabel#1{\strut\vadjust{\kern-\strutdepth\labtekst{#1\ }}}

\def\label #1. #2\par{{\definexref{#1}{\num}{#2}}}
\def\labelf #1\par{{\definexref{#1}{\numfo}{formula}}}
\def\labelse #1\par{{\definexref{#1}{\num}{section}}}


\def\V{\widehat V}
\def\Spec{{\rm Spec}}
\def\A{\widehat{\Bbb A}}
\def\p{{\euscrfm p}}
\def\Z{{\bf Z}}
\def\T{{\bf T}}
\def\X{{\cal X}}
\def\O{{\cal O}}
\def\UA{\overline{A}}
\def\CC{{\Bbb C}}
\def\RR{{\Bbb R}}
\def\QQ{{\Bbb Q}}
\def\PP{{\Bbb P}}
\def\II{{\Bbb I}}

\senu=1 \ssnu=1 \fonu=1

\centerline{\titlefont Analytic subvarieties with many rational
points} \spacing \centerline{ C. Gasbarri}

\bigskip
\bigskip

{\insert\footins{\leftskip\footnotemargin\rightskip\footnotemargin\noindent\eightpoint
$2000$ {\it Mathematics Subject Classification}. Primary 11J97,
14G40, 30D35
\par\noindent {\it Key words}: Transcendence theory, algebraic values of analytic maps, Nevanlinna theory, Bombieri Schneider Lang theorem.}

\vbox{{\leftskip\abstractmargin \rightskip\abstractmargin
\eightpoint

\noindent A{{\sixrm BSTRACT}}.We give a generalization of the
classical Bombieri--Schneider--Lang criterion in transcendence
theory. We give a local notion of $LG$--germ, which is similar to
the notion of $E$-- function and Gevrey condition, and which
generalize (and replace) the condition on derivatives in the theorem
quoted above. Let $K\subset \Bbb C$ be a number field and $X$ a
quasi--projective variety defined over $K$. Let $\gamma\colon M\to
X$ be an holomorphic map of finite order from a parabolic Riemann
surface to $X$ such that the Zariski closure of the image of it is
strictly bigger then one. Suppose that for every $p\in
X(K)\cap\gamma(M)$ the formal germ of $M$ near $P$ is an $LG$--
germ, then we prove that $X(K)\cap\gamma(M)$ is a finite set. Then
we define the notion of conformally parabolic Kh\"aler varieties;
this generalize the notion of parabolic Riemann surface. We show
that on these varieties we can define a value distribution theory.
The complementary of a divisor on a compact Kh\"aler manifold is
conformally parabolic; in particular every quasi projective variety
is. Suppose that $A$ is conformally parabolic  variety of dimension
$m$ over $\Bbb C$ with Kh\"aler form $\omega$ and $\gamma\colon A\to
X$ is an holomorphic map of finite order such that the Zariski
closure of the image is strictly bigger then $m$. Suppose that for
every $p\in X(K)\cap \gamma (A)$, the image of $A$ is an $LG$--germ.
then we prove that  there exists a current $T$ on $A$ of bidegree
$(1,1)$ such that $\int_AT\wedge\omega^{m-1}$ explicitly bounded and
with Lelong number bigger or equal then one on each point in
$\gamma^{-1}(X(K))$. In particular if $A$ is affine
$\gamma^{-1}(X(K))$ is not Zariski dense.

}}

\section Introduction\par

\

A classical theorem by Schneider and Lang (cf. [La]), asserts that
if we have $N>1$ meromorphic functions $f_1(z),\dots, f_N(z)$ which
are of {\it finite order} and such that ${{df_i}\over{dz}}\in
K[f_1,\dots,f_n]$ for some number field $K$, then, if
$Trdeg_{\overline {\Bbb Q}}\overline {\Bbb Q}(f_1\dots, f_N)\geq 2$,
the set of points $z\in\Bbb C$ such that $f_i(z)\in K$ for every
$i$, is finite. This  has been generalized by Bombieri (cf. [Bom]
and [De1]): given $N$ meromorphic functions $f_1,\dots, f_N$ on
$\Bbb C^n$ such that $Trdeg_{\overline{\Bbb Q}}\overline{\Bbb
Q}(f_1,\dots,f_N)\geq n+1$ and with a similar condition on the
derivatives, then the set of points $\underline z\in\Bbb C^n$ such
that $f_i(\underline z)$ is defined and belongs to $K$ for every
$i$, is not Zariski dense in $\Bbb C^N$.

In the Schneider--Lang Theorem, the condition on the derivatives can
be re\-phra\-sed by saying that the image of $\Bbb C$ in $\Bbb C^n$
is the leaf of an algebraic foliation defined over $K$ (and
similarly in the higher dimensional case); on the other side the
fact that the functions are algebraically independent means that
this leaf is Zariski dense in an algebraic variety of dimension
strictly bigger then one.

When one looks closely to the proof of the Schneider--Lang Theorem,
one realize  that the condition on the derivatives, namely that the
ring $K[f_1,\dots, f_n]$ is closed under the derivation, is used
only locally around the points. This condition is needed in order to
derive suitable bounds on the successive derivatives of $P(f_1,\dots
, f_n)$, where $P$ is an "auxiliary polynomial".

Besides, the notion of {\it map of finite order} from a variety $M$
to another variety $X$ make sense in a more general contest then the
classical one appearing in the Bombieri--Schneider--Lang Theorems
(where one deals with meromorphic functions on $\Bbb C^n$).

Actually there is a notion of maps of finite order as soon as we can
define a {it counting function $T_{\gamma}(r)$} \'a la Nevanlinna to
any meromorphic map $\gamma \colon M\to X$. To achieve this it is
natural to assume that $X$ is holomorphically embedded into some
complex projective space, or equivalently (after replacing $X$ by
its Zariski closure) that $X$ is a projective variety. Then, when
$M$ is one dimensional, the counting function may be defined when
the Riemann surface $M$ is {\it parabolic} in the sense of the
classical theory of Myrberg, Nevanlinna, Ahlfors; in particular, as
soon as it is an an affine variety or a finite ramified covering of
the complex plane. Moreover, the extension of Nevanlinna's theory by
Griffiths and King (cf. [GK]) leads to natural notions of
holomorphic maps of finite order from any higher dimensional {\it
smooth affine} complex algebraic variety $M$ to a complex projective
variety $X$. We will see in \S 5 that we can define a value
distribution theory for a wider class of varieties.

It is natural to ask if a general value distribution theory allows
to develop a more general theory of algebraic values of analytic
maps. In this paper we deal with such a generalization.

In the first part of the paper we give a general statement on formal
germs around a rational point (over some number field $K$) on a
projective variety: if some hypothesis are verified then such a germ
is algebraic; these hypothesis are of "Arakelovian" nature: we ask
that the norms (at finite and infinite places) of some maps between
hermitian modules over the ring of integers of $K$ are explicitly
bounded. The benefit of this approach is that, in this way, it is
easier to understand where the difficulties are localized: we can
work on each place independently and then find some global
arithmetic relations.

In the third section we deal with the finite places: we introduce
the notion of $LG$--germ: given a (say, to simplify, one
dimensional) smooth germ $\V$ of analytic variety around a rational
point $P$ of an $N$--dimensional variety $X$, roughly speaking we
say that it is an $LG$--germ of type $\alpha$, if, in suitable
choice of the coordinates around $P$, $\V$  is given by $N$ power
series which are of the form $\sum {{a_iz^i}\over{(i!)^\alpha}}$,
with $a_i$ in the ring of the $S$--integers of $K$, for some finite
set $S$ of places of $K$. The notion of $LG$--germ is exactly what
is needed at finite places in order to let the statement in section
2 work. It is worth remarking that a leaf of a foliation in a smooth
point is an $LG$--germ and to require that a function from $\Bbb
C^n$ to $\Bbb C^N$ defines an $LG$--germ around $K$--rational points
(of $\Bbb C^N$) is less demanding then to require the condition on
the derivatives in the classical Bombieri--Schneider--Lang Theorem:
essentially because the notion of being an $LG$--germ is very local
and the condition on derivatives in the Bombieri--Schneider--Lang
Theorem  is global. The notion of $LG$--germ is similar to the
notion of Gevrey functions and of $E$-- function in the
transcendence theory of values of the solutions of differential
equations developed by Siegel.

In the fourth section we deal with the one dimensional case. Let $K$
be a number field and $\sigma_0:K\hookrightarrow\Bbb C$ is an
embedding of $K$ in $\Bbb C$. We also fix an embedding of the
algebraic closure $\overline K$ of $K$ in $\Bbb C$. Let $X$ be an
$N$ dimensional quasi projective variety defined over $K$; Let
$S\subseteq X(\overline K)$ and, for every positive integer $r$
denote by $S_r$ the set $\{ x\in S\; s. t.\; [\Bbb Q(x):K]\leq r\}$.

The main theorem of the fourth section is:

\label maintheoremforriemannsurfaces1. theorem\par\thm Let $M$ be a
parabolic Riemann surface (with a fixed positive singularity). Let
$\gamma\colon M\to X(\Bbb C)$ be an holomorphic map of finite order
$\rho$ with Zariski dense image. Suppose that, for every $\overline
K$--rational point $P\in S\cap\gamma (M)$, the formal germ $\hat
M_P$, of $M$ near $P$, is (the pull back of) a $LG$-- germ of type
$\alpha$ (in its field of definition). then
$${{Card (\gamma^{-1}(S_r))}\over{r}}\leq {{N+1}\over{N-1}}\rho\alpha[K:\Bbb Q].$$
\endthm

Since every algebraic Riemann surface is parabolic and the leaves of
foliations (defined over $K$) are always $LG$--germs, this give a
generalization of the Schneider--Lang Theorem:

\label slonqbar1. corollary\par\cor Suppose we are in the hypotheses
above, then
$$\sum_{P\in\gamma^{-1}(S)}{{1}\over{[K(\gamma(P));K]}}\leq {{N+1}\over{N-1}}\rho\alpha[K:\Bbb Q].$$
\endcor

If $X$ is a quasi projective variety defined over the number field
$K$ and $r$ is a positive integer, denote by $X_r$ the set $\{ P\in
X(\overline K)\;\; {\rm s.t.}\;\; [\Bbb Q(P):K]\leq r\}$.

\label classicalsl1. corollary\par\cor Let $X$ be an algebraic
variety defined over a number field $K$ and let $F\hookrightarrow
T_X$ be a foliation of rank one (defined over K). Suppose that the
holomorphic foliation $F_{\sigma}\subset (T_X)_{\sigma}$ has a
parabolic leaf $M$ of finite order $\rho$ (for some positive
singularity on $M$) whose Zariski closure has dimension $d>1$, then
$${{Card((X_r\setminus Sing(F))\cap M)}\over{r}}\leq {{d+1}\over{d-1}}\rho[K:\Bbb Q].$$
\endcor

Once one understand the one dimensional case, one realizes that the
main tool used at infinite places is the {\it Evans Kernel}: cf.
definition \ref{evanskernel}. This is a global solution of an
elliptic differential equation and allows to bound the norm of the
jet of a global section of a line bundle $L$ on a point, in term of
the growth of the first Chern form of $L$. This bound is exactly
what is needed for arithmetic applications. Consequently, in the
higher dimensional case one can develop a interesting (for
arithmetic applications) value distribution theory once one suppose
the existence of a function which generalize the properties of the
Evans Kernel.

In the fifth section of the paper we define the notion of {\it
conformally parabolic Kh\"aler ma\-ni\-fo\-lds}. These are Kh\"aler
manifold where one can define an Evans kernel, cf. definition
\ref{laplaceparab}. One should notice the similarity between
parabolic Riemann surfaces and conformally parabolic varieties. We
then show that the main example of conformally parabolic varieties
is the complementary of a divisor in a compact Kh\"aler manifold,
thus in particular every quasi projective variety is conformally
parabolic. We show that we can develop a value distribution theory
for analytic maps from conformally parabolic varieties to projective
varieties. Moreover, again the existence of the Evans kernel allows
to bound the norm of the jet of a global section of a line bundle on
a point, in terms of the growth of the first Chern class of the
bundle. Thus we can develop the arithmetic consequences of this:

If $X$ is a quasi projective variety of dimension $N$ defined over
the number field $K$ and $r$ is a positive integer, denote by $X_r$
the set $\{ P\in X(\overline K)\;\; {\rm s.t.}\;\; [\Bbb Q(P):K]\leq
r\}$. In section 5 we prove:

\label bombierisl1. theorem\par\thm Let $A$ be a $d$ ($d<N$)
dimensional conformally parabolic Kh\"aler manifold with Kh\"aler
form $\omega$. Let $\gamma\colon A\to X(\Bbb C)$ be an analytic map
of finite order $\rho$. Suppose that for every $p\in
X_r\cap\gamma(A)$ the germ of $\gamma(A)$ near $p$ is isomorphic to
an $LG$--germ of type $\alpha$. Suppose that the image of $A$ is
Zariski dense. Then, for every positive integer $r$,  there exists a
current $T_r$ of bidegree $(1,1)$ on $A$ with the following
properties:

--  ${{1}\over{r}}\cdot\int_AT_r\wedge\omega^{d-1}\leq
{{N+1}\over{N-d}}\rho\alpha[K:\QQ]$;

-- For every $p\in\gamma^{-1}(X_r)$, the Lelong number $\nu(T_r,p)$
verifies $\nu(T_r,P)\geq 1$.
\endthm

In analogy with corollary \ref{slonqbar1}, by taking weak limits,
one obtain

\label bslonqbar2. theorem\par\thm  Suppose that we are in the
hypotheses as above. Then there exists a closed positive current $T$
of bidegree $(1,1)$ over $A$ such that:

-- $\int_A T\wedge \omega^{d-1}\leq
2{{N+1}\over{N-d}}\rho\alpha[K:\QQ]$;

-- for every $P\in\gamma^{-1}(X(\overline K)) $ we have that
$\nu(T;P)\geq{{1}\over{[K(\gamma(P)):K]}}$.
\endthm

\endthm

When $A$ is an affine variety, there exists an analogue of
Nevanlinna theory, developed  for instance in  [GK]. Thus there
exists a notion of maps of finite order as soon as we fix an ample
line bundle on a projective compactification $\UA$ of $A$. This
theory  is very similar to the one dimensional case and it is more
or less based on the theory of the Monge--Ampere differential
equations. Unfortunately, at the moment it do not exists an analogue
of the Evans kernel in the theory of the Monge--Ampere equations,
consequently we do not know how to deal with the Cauchy inequalities
over arbitrary affine variety when we are working with "classical"
Nevanlinna theory. For this reason the value distribution theory we
developped is more adapted to the arithmetic problems. In any case
we prove that the order of growth in this theory is the same as the
order computed in the classical theory of Griffiths and King.

As a consequence of \ref{bombierisl1} we find (with the same
hypotheses as before on $X$)

\label bombierischneiderlang1. theorem \par\thm Let $A$ be an affine
variety of dimension $d$ ($d<N$) defined over $\Bbb C$ and
$\gamma\colon A\to X(\Bbb C)$ be an analytic map of finite order
$\rho$. Suppose that for every $p\in S\cap\gamma(A)$ the germ of
$\gamma(A)$ near $p$ is isomorphic to an $LG$--germ of type
$\alpha$. Suppose that the image of $A$ is Zariski dense. Then for
every positive integer $r$, the set $\gamma^{-1}(S_r)$ is contained
in a hypersurface of degree at most
${{N+1}\over{N-d}}{{d!\rho\alpha}\over{\deg(\UA)}}[K:\QQ]r$.\endthm

Once again we can deduce from this an explicit generalization of
Bombieri--Schneider--Lang Theorem (cf. Corollary
\ref{classicalbsl}).

A posteriori one realizes that section 4 is a particular case of
section 5 (thus one may wonder why we decided to write it). One
should notice the following facts: section 4 lean on the preexhising
theory of Parabolic Riemann surfaces and, as written it is
completely self contained, thus a reader who is interested only on
this easier case can skip the more complicated and involved section
5. Methods and definitions of section 5 are inspired by the
corresponding of section 4; we think that the explicit proof in the
one dimensional case allows to better understand and appreciate the
higher dimensional situation. If one want to understand section 5 in
the one dimensional case (and relate it with the classical theory of
parabolic Riemann surfaces), one has to rewrite essentially all
section 4! perhaps up to \ref{schwartzforparabolic} and almost all
what follows. For these reasons we decided to explicitely describe
theory for parabolic Riemann surfaces before we develop the general
case.

This paper owe a lot to J. B. Bost; the section on parabolic Riemann
surfaces is due to him: he kindly explained it to me and sent the
mail [Bo2] where the main ideas and techniques of \S 4 were
explained. In \S 4 we reproduced his ideas. He also spent a lot of
time explaining to me many techniques and ideas: also the
application to the theory of vector bundles with integrable
connections is suggested by him. I warmly thank him for all his
help,  patience and generosity. I thank the referee for the precious
suggestions which highly improve the presentation of  the paper: in
particular the definition of conformally parabolic varieties have
been suggested by him.

\ssection Applications\par

We give an application of Theorem 3.8 based on the recent paper on
classification of foliations on projective surfaces by M. McQuillan
[MQ]: Let $FCP(\Bbb C)$ be the category where the objects are
couples $(Y,f)$ where $Y$ is a Riemann surface and $f\colon Y\to
\Bbb C$ is a finite covering having ramification with polynomial
grown: by this we mean that the function
$$N_f(r):=\sum_{0<\vert f(z)\vert <r}Ord_z(R_f)\log\left|
{{r}\over{f(z)}}\right|$$ where $R_f$ is the ramification divisor,
is a function of polynomial grown ($N_f(r):=O(r^{\rho})$ for some
$\rho$).

\label mcquillan. theorem\par\thm Let $(X,F)$ be a foliated
projective surface defined over a number field $K$ (embedded in
$\CC$). Suppose that the foliation {\rm do not} arises from a
projective connection with irregular singular points and infinite
monodromy. Suppose that $Y$ is a curve in $FCP(\CC)$ contained in
$X$ and invariant for the foliation. Then, if $Y\cap X(K)$ is
infinite, the Zariski closure of $Y$ in $X$ is a curve of genus less
or equal then one.
\endthm

\Proof We may suppose that the Zariski closure of $Y$ is not an
algebraic curve of genus less or equal then one. If the Zariski
closure of $Y$ is an algebraic curve of genus bigger or equal then
two, by the classical Theorem of Faltings (Mordell conjecture) it
has finitely many rational points. If the Zariski closure of $Y$ is
$X$, then, we may look to the classification of 2 dimensional
foliations:  in particular [MQ] chapter V and Fact V.1.2, the
inclusion $Y\to X$ is an analytic map of finite order; so we may
apply theorem \ref{maintheoremforriemannsurfaces1} and conclude.

\smallskip

As an application of Theorem \ref{bombierischneiderlang1} to the
theory of differential equations we can find the following; we
suppose that $K$ is a number field embedded in $\Bbb C$:

\label bost2. theorem\par\thm Let $\overline{X}_K$ be a smooth
projective variety defined over $K$ and $D$ be a simple normal
crossing divisor on it. Denote by $X_K$ the quasi projective variety
$\overline{X}_K\setminus D$. Let $(E;\nabla)$ be a vector bundle on
$\overline{X}_K$ equipped with an integrable algebraic connection
$\nabla$ over $X$ with meromorphic singularities along $D$. Let
$f_{\Bbb C}:(X_K)_{\Bbb C}(\Bbb C)\to E(\Bbb C)$ be an {\rm
analytic} horizontal section of $(E;\nabla)$ defined over $X_{\Bbb
C}:=X_K\times_K\Bbb C$. Let $Y$ be the Zariski closure in $X_K$ of
the set $f^{-1}(E(K))$. Then the restriction of $f_{\Bbb C}$ to $Y$
is an {\rm algebraic} horizontal section of $(E,\nabla)\vert_Y$.
\endthm

\Proof We may suppose that $X_K$ is smooth and affine. Taking as
$X_K$ a desingularization of an irreducible component of $Y$, we may
suppose that $f^{-1}(E(K))$ is Zariski dense and eventually prove
that $f_{\Bbb C}$ is algebraic. Remark that $f_{\Bbb
C}^{-1}(E(K))\subseteq X_K(K)$, consequently to require that
$f_{\Bbb C}^{-1}(E(K))$ is Zariski dense on $X_K$ is equivalent to
require that it is Zariski dense on $X_K\times_K\Bbb C$. We claim
that the section $f_{\Bbb C}$ is of finite order:

We may work in a neighborhood $U$ of  a point of the divisor $D$. We
may suppose that $U$ is a polydisk with coordinates $z_1,\dots, z_n$
and that $D$ is given by $z_1\cdot z_2\cdot\dots \cdot z_r=0$.

The connection is given by the linear differential equation
$$dy={{A}\over{z_1^{i_1}\cdot\dots z_r^{i_r}}}y$$
With $A$ a matrix of holomorphic forms.

Let $z=(z_1,\dots, z_n)$ with $\vert z_i\vert=1$ and $t\in [0,1]$
and consider the function $h(t):=\vert y(tz)\vert$. By Lipschitz
condition, there is a suitable constant $C$ such that

$$\vert h'(t)\vert\leq \vert z\cdot dy(tz)\vert=\vert dy(tz)\vert\leq {{C}\over{t^{i_1+\cdots + i_r}}}h
(t).$$

Observe that $h(t)$ is decreasing, thus, $h'(t)<0$. Diving by $h(t)$
and integrating both sides we obtain
$$\log(h(t))\leq {{C'}\over{t^a}}$$
for suitable $C'$ and $a$. The claim follows.

We apply now Theorem \ref {bombierischneiderlang1} to $X_K$ and we
obtain that, if $f_{\Bbb C}$ is not algebraic, $f^{-1}(E(K))$ cannot
be Zariski dense. From this we conclude.

\smallskip

\rmk One should notice that, an analytic function, say over $\CC^d$,
is of finite order $\rho$ if there exists a constant $C>0$ such that
$\sup_{\Vert z\Vert\leq R}\{ \vert f(z)\vert\}\leq C^{R^\rho}$. This
definition is not the same as the definition given in [Dl].

\endssection

\

\endsection

\

\section The general algebraization setting\par

\

Let $K$ be a number field, $O_K$ be its ring of integers,
$S_{\infty}$  the set of infinite places, $M_{fin}$ the set of
finite places of $K$ and $M_K=S_{\infty}\cup M_{fin}$. We fix once
for all an embedding $\sigma_0 :K\hookrightarrow \Bbb C$ (so an
infinite place of $K$). Let $X_K$ be a geometrically irreducible,
connected quasi projective variety of dimension $N$ defined over $K$
and $F\subseteq X(K)$ a finite subset. For $P\in F$ we denote by
$\widehat{X}_P$  the completion of $X_K$ at $P$. for every $P\in F$,
let $\V_P\subseteq \widehat{X}_P$ be a smooth formal subvariety of
dimension $d$. For every (finite or infinite) place $\euscrfm p$ of
$K$, we suppose that the base change $V_P^{\euscrfm p}$ is analytic:
the power series defining it have  radius of convergence which is
non zero.

\label algebraicformalscheme. definition\par\defi Let $P\in X(K)$
and $\V_P\subseteq\widehat{X}_P$ be a formal subscheme. We will
denote by $\bar{V_P}$ the smallest {\rm Zariski} closed set, defined
over $K$, containing $\V_P$ and call it {\rm the Zariski closure of
$\V_P$}.

We will say that $\V_P$ {\rm is algebraic} if
$\dim\V_P=\dim\bar{V_P}$.
\enddefi

In this chapter we would like to give some general sufficient
condition in order to have that each of the $\V_P$'s is algebraic.

By replacing $X_K$ with the Zariski closure of the $\V_P$'s we may
suppose that at least one of the $\V_P$'s is Zariski dense.

We fix a model $\X$ of $X_K$ flat, projective over $\Spec (O_K)$.
Moreover for every infinite place $\sigma$ we fix a K\"ahler metric
on $X_{\sigma}$.

Let $L$ be a relatively ample line bundle on $\X$ and, we suppose
that, for every infinite place $\sigma$ the holomorphic line bundle
$L_{\sigma}$ is equipped with a smooth positive metric.

Consequently, for every positive integer $D$, the locally free
$O_K$--module $E_D:=H^0(\X; L^D)$ is equipped with the $L^2$ and the
$\sup$ norms which are comparable by, for instance [Bo]; so $E_D$
has the structure of an hermitian vector bundle over $O_K$.

For every $P\in F$ we fix an hermitian integral structure on the
$K$--vector space $T\V_P$. Consequently, for every place $\euscrfm
p$ of $K$, and for every $(i,D)\in\Bbb N\times\Bbb N$, the vector
space $\left(S^i(T^\ast_P\V)\otimes L^D\vert_P)\right)_{\euscrfm p}$
is equipped with a $\euscrfm p$--adic norm  (which, at infinite
places is a metric).

Let $(V_{P})_i$  the $i$--th infinitesimal neighborhood of $P$ in
$\V_{P}$; there is a canonical exact sequence

$$0\longrightarrow S^i(T^\ast_{P}\V)\otimes L^D\vert_{P})\longrightarrow
L^D\vert_{(V_{P})_{i+1}}\longrightarrow L^D\vert_{(V_{P})_{i}}.$$

For every $(i;D)\in \Bbb N\times \Bbb N$, we denote by $E^i_D$ the
kernel of the map $\alpha^i_D$ given by the composite \advance\ssnu
by-1\labelf mapalpha\par$$E_D\hookrightarrow \bigoplus_{P_j\in
F}H^0(\V_{P_j};L^D)\longrightarrow\bigoplus_{P_j\in
F}H^0((\V_{P_j})_i;L^D)\eqno{{(\numfo)}}$$\advance\ssnu by1

remark that the first map is injective because the germs are Zariski
dense. We denote by $\gamma^i_D$ the induced map \labelf
mapgamma\par$$\gamma^i_D:E^i_D/E^{i+1}_D\hookrightarrow
\bigoplus_{P_j\in F}H^0(P_j; S^i((T_{P_j}\V_{P_j})^{\ast}) \otimes
L^D\vert_{P_j}).\eqno{{(\numfo)}}$$\advance\ssnu by1

For every (finite or infinite) place $\euscrfm p$, the $K_{\euscrfm
p}$--vector space $(E^i_D)_{\euscrfm p}$ is equipped with the norm
induced by the norm of $(E_D)_{\euscrfm p}$. Consequently the
$K_{\euscrfm p}$--vector space
$\left(E^i_D/E^{i+1}_D\right)_{\euscrfm p}$ is equipped with the
quotient norm.  We will denote by $\Vert\gamma_D^i\Vert_{\euscrfm
p}$ the  $\euscrfm p$--norm of the linear map $\gamma_D^i$.

An easy computation show that ${{1}\over{[K:\Bbb Q]}}\sum_{\sigma\in
M_K}\log\Vert \gamma\Vert_{\sigma}$ do not depend on the field $K$.

\label generalalgebraization. theorem\par\thm With the notation as
above, suppose that the formal germs are {\rm not algebraic} and
that we can choose the hermitian integral structure on the tangent
spaces $T_P\V$'s in such a way that the following holds: there exist
constants $C_1$, $C_2$, $\lambda$ and $A$ depending only on $F$, the
model $\X$ and the choice of the hermitian integral structure, but
independent on the $i$ and $D$, for which:

1) the following inequality holds: $${{1}\over{[K:\Bbb
Q]}}\sum_{\euscrfm p\in M_{fin}}\log\Vert \gamma_D^i\Vert_{\euscrfm
p}\leq C_1(i\log i+ C_2(i+D));$$

2) the following inequality holds
$${{1}\over{[K:\Bbb Q]}}\sum_{\sigma\in S_{\infty}}\log\Vert\gamma_D^i\Vert_{\sigma}\leq C_2(i+D);$$

3) for ${{i}\over{D}}\geq \lambda$, the following inequality holds€
$${{1}\over{[K:\Bbb Q]}}\sum_{\sigma\in S_{\infty}}\log\Vert\gamma_D^i\Vert_{\sigma}\leq
C_1\left(-Ai\log\left({{i}\over{D}}\right)+ C_2(i+D)\right).$$

Then $A\leq{{N+1}\over{N-d}}$ (recall that $N= dim(X)$ and
$d=\dim(V_P)$).
\endthm

At first glance theorem \ref{generalalgebraization} may seems
strange. One should understand it in this way: given some formal
germs, if one can prove that the norms of the involved maps are {\it
very small} then they are algebraic. One will see in the other
sections that the bounds of the norms is related to the arithmetic
and the geometry of the formal germs.

Before we start the proof, we recall the main slope inequalities,
proved, for instance in [Bo], \S 4.1 (we refer to loc. cit. for the
notation):

a) If $E$ is an hermitian vector bundle over $O_K$, then we call the
real number $\mu_n(E):={{1}\over{[K;\Bbb
Q]}}\cdot{{\widehat{\deg}(E)}\over{rk(E)}},$ {\it the slope} of $E$;

b) within all the sub bundles of a given hermitian vector bundle
$E$, there is one having maximal slope; we call its slope {\it the
maximal slope of $E$} and denote it by $\mu_{\max}(E)$;

c) If $E$ is an hermitian vector bundle and $S^i(E)$ its $i$--th
symmetric power metric. we endow  $S^i(E)$ with the symmetric power
metric of $E$.  We can find a constant (depending on $E$) $B\in\Bbb
R$, such that $\mu_{\max}(S^i(E))\leq iB$;

d) if $E_1$ and $E_2$ are two hermitian vector bundles, then
$$\mu_{\max}(E_1\oplus E_2) = \max \{\mu_{\max}(E_1);
\mu_{\max}(E_2)\};$$

e) (the main slope inequality): Suppose that $F_K$ is a finite
dimensional vector space endowed with a filtration $\{
0=F^{r+1}\}\subset F^r\subset F^{r-1}\subset\dots\subset F^0=F_K$
such that the successive sub quotients $G^i_K=F^i/F^{i+1}$ are the
generic fibre of hermitian vector bundles $G^i$. Suppose moreover
that we have an hermitian vector bundle $E$ and an {\it injective}
linear map $\varphi_K\colon E_K\hookrightarrow F_K$. Denote by
$E^i_K:=\varphi^{-1}(F_i)$ and by $E^i:=E^i_K\cap E$. For every $i$,
we have an induced map $\varphi_i\colon E^i\to G^i$. For every
(finite or infinite) place $\euscrfm p$, we denote by $\Vert
\varphi_i\Vert_{\euscrfm p}$ the  $\euscrfm p$--adic norm  of the
linear map $\varphi_i$; then
$${{\widehat{\deg}(E)}\over{[K;\Bbb Q]}}\leq \sum_{i=0}^rrk\left(
E^i/E^{i+1}\right)\left[\mu_{\max}(G^i)+\sum_{\euscrfm
p}\log\Vert\varphi_i\Vert_{\euscrfm p}\right];$$

f) (Arithmetic Hilbert--Samuel formula) There exists a positive
constant $C$ such that
$${{\widehat{\deg}(E_D)}\over{[K;\Bbb Q]}}\geq -CD^{N+1}$$
(where we recall that $N=\dim X_K$).

We can now prove theorem \ref{generalalgebraization}:

\Proof {\it (of Theorem \ref{generalalgebraization})} We suppose
that at least one of the $\V_P$'s  is {\it not} algebraic and that
$A>{{N+1}\over{N-d}}$ and we eventually find a contradiction.

Choose an $\epsilon >0$, such that
$\alpha:={{N+1}\over{d+1}}-\epsilon>1$ and $\alpha>{{A}\over{A-1}}$.

Take $D$ such that $D^{\alpha}>\lambda D$. Then, by the main slope
inequality applied to the map \ref{mapalpha} with the filtration
induced by \ref{mapgamma} and the arithmetic Hilbert--Samuel
formula, we can find constants $C_j$'s, independent on $i$ and $D$,
(provide that they are sufficiently big) such that
$$\eqalign{-C_0D^{N+1}&\leq\sum_{i=0}^{\infty}rk\left(
E_D^i/E_D^{i+1}\right)\left[C_2(i+D)+i\log i+\sum_{\sigma\in
S_{\infty}}\log\Vert\gamma_D^i\Vert_{\sigma}\right]\cr
&\leq\sum_{i\leq D^{\alpha}}rk \left(
E_D^i/E_D^{i+1}\right)\left[C_3D^{\alpha}+C_4D^{\alpha}\log
D\right]\cr & +\sum_{i> D^{\alpha}}rk\left( E_D^i/E_D^{i+1}\right)
\left[ C_5\left(-A i\log\left( {{i}\over{D}}\right)+i\log
i\right)+C_2(i+D)\right].\cr}$$

Consequently
$$\eqalign{C_0&\geq {{rk\left( E_D/E_D^{D^{\alpha}}\right)
\left[-C_3D^{\alpha}-C_4D^{\alpha}\log D\right]}\over{D^{N+1}}}\cr &
+ \sum_{i>D^{\alpha}}{{rk\left(
E_D^i/E_D^{i+1}\right)}\over{D^N}}\left(C_5{{i}\over{D}}\left[(A
-1)\log i-A\log
D\right]-C_2\left({{i}\over{D}}+1\right)\right).\cr}$$

If we denote by $I_j$ the ideal of the point $P_j$ on the local germ
$\V_{P_j}$, we have an injection
$$E_D/E_D^{D^{\alpha}}\hookrightarrow \bigoplus_{P_j\in
F}L^D\vert_{P_j}\otimes\O_{\V_{P_j}}/I_j^{D^{\alpha}};$$

thus $rk(E_D/E_D^{D^{\alpha}})=O(D^{d\alpha})$. This, together with
our choice of $\alpha$, implies that
$$\lim_{D\to\infty}{{rk(E_D/E_D^{D^{\alpha}})
\left(-C_3D^{\alpha}-C_4D^{\alpha}\log D\right)}\over{D^{N+1}}}=0.$$

On the other side, for $i\geq D^{\alpha}$, we have
$$(A -1)\log i\geq \alpha(A-1)\log D;$$
consequently, again because of our choice of $\alpha$, we can find a
positive constant $\beta$ such that
$$(A -1)\log i-A \log D>\beta\log D.$$
Thus
$$\eqalign{&\sum_{i>D^{\alpha}}{{rk\left(
E_D^i/E_D^{i+1}\right)}\over{D^N}}\left(C_5{{i}\over{D}}\left[(A
-1)\log i-A\log D\right]-C_2\left({{i}\over{D}}+1\right)\right)\cr &
\sum_{i>D^{\alpha}}{{rk\left(
E_D^i/E_D^{i+1}\right)}\over{D^N}}\left[C_5{{i}\over{D}}\beta\log
D-C_2\left({{i}\over{D}}+1\right)\right]\cr &\geq
\sum_{i>D^{\alpha}} {{rk\left(
E_D^i/E_D^{i+1}\right)}\over{D^N}}\left[C_5{{i}\over{D}}\left(\beta\log
D-C_6\right)-C_2\right].\cr}$$

If we choose $D$ so big such that $\beta\log D-C_6\geq 0$ we obtain
 that the last displayed sum is greater than
$$\eqalign{&C_6\sum_{i>D^{\alpha}}{{rk\left(
E_D^i/E_D^{i+1}\right)}\over{D^N}}\left[ D^{\alpha-1}-1\right]\cr &
=C_6{{rk(E_D^{D^{\alpha}})}\over{D^N}}\left(
D^{\alpha-1}-1\right);\cr}$$

but since $\lim_{D\to\infty}{{rk(E_D^{D^{\alpha}})}\over{D^N}}=1$,
we eventually find a contradiction as soon as $D$ is sufficiently
big.

\smallskip

Once we assure the conditions (1), (2) (3), Theorem
\ref{generalalgebraization} gives us a powerful tool for
algebraization of analytic germs of subvarieties of a projective
variety. We remark that the classical Cauchy inequality assures
that, if the formal germ has a positive radius of convergence at
every infinite place, condition (2) is automatic. On the other side
conditions (1) and (3) not always hold. If one chose the integral
structures involved in a an arbitrary way there is no way a priori
to give explicit bounds as required by the hypotheses of the
theorem. One may think the problem in this way: the choices of the
metrics at infinity may be done by some analytic information of the
germs; the choice of the integral structure is related to $p$--adic
information of the germs and the choice of the constant $A$ is
purely arithmetic.

\endsection

\

\section $LG$--germs\par

\

We describe in this chapter a class of germs of analytic
subvarieties of a variety defined over a number field which verify
condition (1) of Theorem \ref{generalalgebraization}. Roughly
speaking germs in this class are defined by power series whose
coefficients of each monomial are algebraic integers divided by the
factorial of its exponents. the example to keep in mind is the
classical exponential function. The other leading example is the
formal leaf of an analytic foliation in a smooth point of the
foliation.

At the moment we do not know the right condition a germ should have
in order to apply theorem \ref{generalalgebraization} (beside the
fact that it must verify (1) of loc. cit.). The referee proposed
(and we warmly thank her/him for this) very interesting remarks and
observations on the Gevrey condition on analytic power series. We
hope to come again on this in a future paper.
\smallskip

\label goodgerms. section\par\ssection Good germs\par Let $K$ be our
number field and $S\in M_K$ be a finite set of places of $K$
(containing all the infinite places); let $O_S$ be the ring of
$S$--integers of $K$:
$$O_S:=\left\{ a\in K/\,\,\Vert a\Vert_v\leq 1\,{\rm for\; every
} v\not\in S\right\}.$$

We will put $B:=Spec (O_S)$. In general if $R$ is scheme, we will
denote by $\Bbb A^N_R$ the $N$-- dimensional affine space over $R$
and by $\A^N_R$ the completion of it at the origin $0$.

Let $f:\X\to B$ be an irreducible scheme, flat over $B$ and ${\cal
P}:B\to\X$ be a rational point. We will suppose that $f$ is smooth
in a neighborhood of $ {\cal P}(B)$. When $f_K:{\X}_K\to Spec (K)$
is smooth around ${\cal P}_K$, this will be the case up to enlarge
$S$, if necessary.

Let $N=\dim ({\X}_K)$.

Let $\widehat{{\X}_P}$ be the completion of $\X_K$ around $P:={\cal
P}_K$.

Since $\X$ is smooth over $B$ around ${\cal P}$, we can choose an
open neighborhood $U\subset\X$ of $P$ and an \'etale map $g:U\to
\Bbb A_B^N$. We can (and we will) suppose that $g({\cal P})=0$.

By definition of \'etale map, the completion of $g$ will induce an
isomorphism $\hat g:\widehat{\X_P}\buildrel\sim\over\to
\widehat{\Bbb A_B^N}$.

We fix coordinates ${\Z}=(Z_1;\dots ;Z_N)$ on $\A^N_B$ and
$\T:=(T_1;\dots;T_d)$ on $\A^d_B$. This choice induces coordinates
on $\A^N_K$ and on $\A^d_K$.

Let $\gamma: \A^d_K\to\widehat{\X_P}_K$ be a formal morphism (such
that $\gamma (0)=P$). The map $h:=\hat g\circ \gamma\colon
\A^d_K\to\A^N_K$ is given by $N$ power series
$$h(\T)=\Z(\T)=(Z_1(\T);\dots ; Z_N(\T))$$
with $Z_i(\T):=\sum_{I\in \Bbb N^d}a_I(i){\T}^I$ and $\Z(0)=0$.
Where, here and in the sequel, if $I=(i_1,\cdot , i_d)\in \Bbb N^d$,
we will denote $\T^I:=T_i^{i_1}T_2^{i_2}\cdots T_d^{i_d}$ moreover
we will write $I!:=i_1!\cdot i_2!\cdot\dots\cdot i_d!$ and $\vert
I\vert:=i_1+\cdot +i_d$.

\smallskip

\label defofsgoodgerm. definition\par\defi Let $\alpha$ be a real
number. The map $\gamma$ will be said to be a {\rm $S$--Good germ of
type $\alpha$} with respect to $\X$ (or simply a $S$--Good germ, if
$\alpha$ is clear from the contest), if:

a) For every $\p\not\in S$, there in a constant $C_\p\geq 1$ such
that
$$\Vert a_I(i)\Vert_\p\leq {{C_\p^{\vert I\vert}}\over{\Vert I!\Vert_\p^{\alpha}}} \;\;\;\;\; {\rm and }\;\;\;\; \prod_{\p\not\in
S }C_\p<\infty.$$

b) The $a_I(i)$ with $\vert I\vert=1$ are in $O_S$ and the linear
map $\alpha:O_S^d\to O_S^N$ given by the matrix $(a_I(i))$ with
$\vert I\vert=1$ is is injective and has cokernel without torsion.
\enddefi

The second condition can be geometrically explained in the following
way: we have an induced linear map of tangent spaces
$$d\gamma :T_0(\widehat {\Bbb A_B^d}))\longrightarrow T_P\X_K$$ in both tangent spaces there is a natural $O_S$ module of
maximal rank given by the tangent space of their integral models.
Condition  b) is equivalent to ask that the map $d\gamma$ is the
extension of a linear map between these $O_S$--modules which is
injective modulo every prime ideal of $O_S$. We observe that a
$S$--good germ is a smooth formal subscheme of $\widehat{\X_{P_K}}$.

In the sequel we will use the following lemma

\lemma For every positive integer $i$, the following holds
$$-\sum_{\p\in M_{fin}}\log\Vert i!\Vert_{\p}\leq i\log(i)-i+O(1).$$\endlemma

The proof is an easy consequence of the product formula and the
Stirling's formula (we leave details to the reader).

\smallskip

The power series $Z_i(\T)$ define a $S$--Good germ if and only, for
every $\p$ there is a constant $C_\p$ with the same properties as
definition \ref{defofsgoodgerm}, such that, for every multiindex $I$
(with the evident notation)
$$\Vert{{\partial^{I}Z_i(T)}\over{\partial
T^I}}(0)\Vert_\p\leq {{C_\p^{\vert I\vert}}\over{\Vert
I!\Vert^{\alpha-1}_\p}}.$$

A useful property of $S$--Good germs is deduced from the following
lemma, the proof of which is done by easy induction:

\label productofgoods. lemma\par\lemma Let $W_i:=\sum_Ia(i)T^I\in
K[\![\T]\!]$ ($i=1,2$) be two power series. Suppose that there are
constants $C_i$ such that
$$\Vert a_I(i)\Vert_\p\leq {{C_i^{\vert I\vert}}\over{\Vert I!\Vert_\p^{\alpha}}}.$$ Then  $W_1\cdot W_2=\sum_Ib_IT^I$ with
$\Vert b_I\Vert_\p\leq {{(C_1\cdot C_2)^{\vert I\vert}}\over{\Vert
I!\Vert_\p^{\alpha}}}.$
\endlemma

The property of being $S$-- good, for a germ $\gamma$, seems to
depend on the choice of the coordinates on $\A^d_S$, on $\A^N_S$ and
on the map $g$; so it would be more precise to speak about
"$S$--goodness with respect to $(\T, g,\Z)$". As the lemmas below
will show, the fact that a map $\gamma$ is a $S$--good germ depends
only on $\gamma$ itself and on the model $\X$. So we can speak about
$S$-- good germs over a flat $S$--scheme without making reference on
any choice.

\lemma Let
$$\eqalign{Aut_S(\A_B^d):=&\{ F(\T)=(f_1(\T);\dots ;f_d(\T))\,  s.t.\cr & f_i(\T)=\sum_{J\in\Bbb
N^d}b_JT^J, \, \, b_J\in O_S\,\, and\,\, \det({{\partial
f_i}\over{\partial t_j}})(0)\in O_S^\ast\}\cr}$$

be the group of $S$--automorphisms of $\widehat {\Bbb A_B^d}$. If
$\T=F(\T')$ for some $F\in Aut_S(\widehat {\Bbb A_S^d})$, then the
germ $\gamma$ is $S$--good with respect to $(\T; g; \Z)$ if and only
if it is $S$--good with respect to $(\T'; g; \Z)$.\endlemma

The definition of $Aut_S(\widehat {\Bbb A_S^d})$ is intrinsic it do
not depend on the choice of the coordinates: for every finite place
$\p$ of $O_S$ denoting by $\Delta^d_{\p}$ the poly--disk $\{
(x_1;\dots ;x_d)\in K_{\p}^d\, /\, \Vert x_i\Vert_{\p}< 1\}$;
$Aut_S(\widehat{\Bbb A^d_S}))$ the is the group of the formal
automorphisms $F$ of $\widehat {\Bbb A_K^d}$  such that, for every
$\p\not\in S$ the restriction of $F$  to $\A^d_{K_\p}$ converges on
$\Delta^d_{\p}$ and $F(\Delta^d_{\p})=\Delta^d_{\p}$.

\Proof With the notations as in \ref{defofsgoodgerm}, given a good
germ $\gamma$ and an automorphism $F\in Aut_S(\A^d_S)$ we have a
composite of maps
$$\widehat {\Bbb A_S^d}\buildrel F\over\longrightarrow \widehat {\Bbb A_S^d}
\buildrel h\over\longrightarrow \widehat {\Bbb A_S^N}$$ and the
corresponding map between tangent spaces. Since $F\in Aut_S(\widehat
{\Bbb A_S^d})$, the linear map $dF$ is an isomorphism of the
integral models, so condition b) holds for $T'$ (and for $g$ and
$\Z$).

We show now that, if condition (a) holds for $(\T; g;\Z)$ it holds
for $(\T'; g;\Z)$. Let $\Z_i(\T')=\Z_i(F(\T'))$ we must show that,
there exists  constants $C_\p$ (with the condition on all the
$C_\p$'s) and $\alpha$ such that, for every $I=(i_1;\dots ,
i_d)\in\Bbb N^d$,
$$\Vert{{\partial^I}\over{(\partial^I T')}}\Z(0)\Vert_\p\leq {{C_\p^{\vert I\vert}}\over{\Vert I!\Vert_\p^{\alpha}}}.$$

by induction on the multi--index $I$ (ordered by lexicographic
order) it is easy to see that for a given $I\in \Bbb N^d$ there
exists a polynomial $P(X_{\alpha};Y^h_{\beta})\in\Bbb
Z[X_{\alpha};Y^h_{\beta}]$  where $\alpha\in \Bbb N^d$ is such that
$\alpha\leq I$ and $\beta$ and $h$ are suitable multi--index having
the following properties:

\noindent i) $\deg P(X_{\alpha}; Y_{\beta}^h)_{X_\alpha}\leq 1$;

\noindent ii)  ${{\partial^I}\over{(\partial^I
T')}}\Z(\T')=P({{\partial^{\alpha}\Z(F(\T'))}\over{(\partial
\T)^{\alpha}}};{{\partial^hF_{\beta}(\T')}\over{(\partial
\T')^h}})$.

If we put $\T'=0$ and use the hypotheses we  easily conclude.

\rmk Observe that, the constants $C_\p$ may depend on the
parametrization but the constant $\alpha$ depends only on the
germ.\endrmk

\smallskip

Using the same strategy we prove that the notion of $S$--Good germ
independent on the choice of the coordinates on $\A^N_S$:

\label independence. lemma\par\lemma  If ${\bf W}=G(\Z)$ for some
$G\in Aut_S(\widehat{A_S^N})$ then the germ $\gamma$ is $S$--Good
with respect to $(\T; g;{\bf W})$ if and only if it is $S$--Good
with respect to $(\T; g; \Z)$.\endlemma

The fact that the notion of $S$--good germ is independent on the
choice of $g$ follows also from Lemma \ref{independence}.

\ssection Definition of $LG$--germ\par

Let $X_K$ be a smooth variety defined over $K$ of dimension $N$ and
$p\in X_K(K)$. Let $\widehat {X_P}$ be the formal neighborhood of
$X_K$ around $P$ and $\iota:\V\hookrightarrow \widehat{X_P}$ be a
smooth formal subscheme of dimension $d$.

\label LG-germ. definition\par\defi Let $\alpha$ be a real number.
The formal scheme $\V$ is a {\rm $LG$--germ of type $\alpha$} (or
simply a $LG$--germ) if the following holds:

i) For all $\p\in \Spec\max(O_K)$ the induced formal subscheme
$\V_{K_{\p}}\hookrightarrow (\widehat{X_P})_{K_{\p}}$ is an analytic
germ: the equations defining it have a positive radius of
convergence.

ii) There exists a finite set $S$ of places of $K$ (containing all
the infinite places) with:

-- a smooth model $\X_S\to Spec (O_S)$ of $X_K$ over which the
rational point extends to a section $P:Spec(O_S)\to \X_S$,

-- a $S$--good germ of type $\alpha$ with respect to $\X$,  $\gamma
:\A_K^d\hookrightarrow\widehat{\X_S}_P$

-- a $K$--isomorphism $\delta :\A_K^d\buildrel\sim\over\to\V$ such
that $\gamma=\iota\circ\delta$.
\enddefi

\rmk Suppose that $\iota:\V\to(\widehat{X_K})_P$ is an $LG$--germ
and $\X\to Spec (O_K)$ is a model of $X_K$. then we can find a
finite set of places $S$ a model $\X'$ of $X_K$, a birational map
$\X' \to\X$ and a $S$--good germ with respect to $\X'$,
$\gamma:\A_K^d\hookrightarrow\widehat{X_P}$ with a $K$--isomorphism
$\delta:\A_K^d\buildrel\sim\over\to\V$ such that
$\gamma=\iota\circ\delta$. consequently, the notion of $LG$--germ is
essentially independent on the choice of the model.
\endrmk

One of the leading example of $LG$-germ is the formal leaf of a
foliation in a smooth point:

\label foliationlg. propositon\par\prop $X_K$ be a smooth variety
over $K$ and ${\cal F}\hookrightarrow T_{X_K}$ be a foliation on it.
Let $P\in X_K(K)$ be a smooth point for the foliation and
$V_P\hookrightarrow \widehat{X_K}_P$ the formal leaf of the
foliation through $P$. Then $V_P$ is a $LG$--germ of type 1.\endprop

\prop. Suppose we are in  the hypotheses of \ref{foliationlg} and
moreover for almost every $\p$ the foliation is closed under the
$p$--derivation; then the formal leaf of the foliation through a
smooth point is a $LG$--germ of type zero.\endprop

The proofs are in [Bo] page 189 and ff.

\smallskip

Let $X_K$, $P\in X_K(K)$ and $\V\hookrightarrow\widehat{X_P}$ be a
smooth formal germ. An integral model of $X_K$ give rise to an
integral structure on the tangent space of $\V$:

\label integralstructure. \par\ssection Construction:\par Let
$f:\X\to\Spec (O_K)$ be a model of $X_K$ where $P$ extends to a
section ${\cal P}:Spec (O_K)\to\X$. Denote by $T^\ast(\cdot)$ the
cotangent space. By definition we have a surjective map
$(d\iota)^\ast:T_P^\ast X_K\twoheadrightarrow T_P^\ast\V$  and
$T^\ast_PX_K$ is the generic fibre of the $O_K$-- module
$T^\ast_{\cal P}\X$; the $O_K$-- module $(d\iota)^\ast(T^\ast_{\cal
P}\X)$ is then an integral model of $T^\ast_P\V$. We will denote it
by ${\cal T}_P^\ast\V$.\endssection

Let $X_K$ and $P\in X_K(K)$ as above. Let
$\V\hookrightarrow\widehat{X_P}$ be a $LG$--germ. Denote by $\V_i$
the $i$--th infinitesimal neighborhood of $P$ in $\V$. There is a
natural exact sequence
$$0\longrightarrow S^i(T^\ast_P\V)\longrightarrow\O_{\V_{i+1}}\longrightarrow\O_{\V_{i}}\longrightarrow
0.$$ Let $U\subset X_K$ be a Zariski open neighborhood of $P$ and
$s\in H^0(U, \O_{X_K})$. Suppose that the restriction of $s$ to
$\V_{i-1}$ is zero (we will say that $s$ vanishes at order $i-1$
along $\V$). Then, the restriction of $s$ to $\V_i$ canonically
defines a section in  $S^i(T^\ast_P\V)$. We will denote this section
by $j^i(s)$ and call it the {\it $i$--th jet of $s$}.

Let $\X\to\Spec(O_K)$ be a model of $X_K$ and suppose that $P$
extends to a section $P\colon\Spec(O_K)\to\X$.

The interest of $LG$--germs is that, we can actually bound from
above the norm of the $i$--th jet of integral functions.

\label boundofnorms. theorem\par\thm Let $\V\hookrightarrow
\widehat{(X_K)}_P$ be a $LG$--germ of type $\alpha$ and
$\X\to\Spec(O_K)$ be a model of $X_K$. Then we can find a constant
$C$ for which the following holds:

Let $U\subseteq\X$ be a Zariski open neighborhood of $P$ and $s\in
H^0(U,\O_\X)$. Suppose that $s$ vanishes at the order $i-1$ along
$\V$ then
$$\sum_{\p\in\Spec\max(O_K)}\log\Vert j^i(s)\Vert_\p\leq \alpha
i\log(i)+Ci.$$

\endthm

\Proof First of all, we remark that, since $\V$ is a $LG$--germ we
can work with the involved map $\gamma\colon\A^d_K\to
\widehat{(X_K)}_P$. We fix the set $S$ and we can suppose that the
smooth model $\X'\to B:=\Spec(O_S)$  involved in the definition of
the $LG$--germ is $\X\vert_B$. We fix an open neighborhood
$U\subseteq\X\vert U$ of $P$ with an \'etale map $g\colon U\to {\Bbb
A}^N_B$. Fix coordinates $\Z:=(Z_1,\dots ,Z_N)$ on $\A^N_B$ and
coordinates $\T:=(t_1,\dots ,t_d)$ on $\A^d_B$. The restriction to
$B$ of the integral structure of $S^i(T^\ast_P\V)$ is the
$O_S$--module ${\cal S}^i:=\bigoplus_{i_1+\dots
+i_d=i}O_S(dt_1)^{i_1}\cdots (dt_d)^{i_d}$.

Moreover the restriction of $s$ to $\A^N_B$ is  a power series
$s=\sum_{I\in\Bbb N^N}a_IZ^I$ with $a_I\in O_S$. The restriction via
$\gamma$ of $s$ to $\A^d_K$ is then the power series
$s\vert_{\A^d_K}=\sum a_I\Z^I(\T)$. Thus, considering Taylor
extension,
$$j^i(s)=\sum_{i_1+\dots +i_d=i}{{1}\over{i_1!\cdot i_2!\cdot\dots\cdot
i_d!}} \cdot{{\partial^is(\T)}\over{\partial^{i_1}t_i\dots
\partial^{i_d}t_d}}\vert_{\T=(0)}(dt_1)^{i_1}\cdots
(dt_d)^{i_d}.$$ The conclusion follows by induction on $I$ and
\ref{productofgoods} applied to $Z^I$.

\smallskip

As corollary we find the estimation at finite places needed to apply
Theorem \ref{generalalgebraization}. Let $X_K$ be a smooth
projective   variety defined over $K$, $P\in X_K(K)$ a $K$--rational
point and $\iota:\V\to X_K$ be a $LG$--germ. Suppose that $L$ is an
ample line bundle on $\X$ (a suitable model of $X_K$) equipped with
a positive hermitian metric. We suppose that we fixed a positive
metric on $X_K({\Bbb C})$. Let $\gamma_D^i$ the linear maps defined
in \ref{mapgamma}; then

\label finitenorms. corollary\par\cor Let $X_K$ be a smooth variety
defined over $K$, $P\in X_K(K)$ and $L$ as above. Let
$\V\subseteq\widehat{X_P}$ be a smooth formal subvariety. Suppose
that $\V$ is an $LG$--germ of type $\alpha$. Then for every model
$\X$ of $X_K$ over $O_K$ we can choose an integral model of the
tangent space of $\V$ at $P$ and a positive constant $C$ depending
only on the model, on the point $P$ and on $\V$ such that
$$\sum_{\p\in M_{fin}}\log\Vert\gamma_D^i\Vert_{\p}\leq
[K:Q]\left(\alpha\cdot i\log i+Ci\right).$$
\endcor
In order to prove this corollary we use Theorem \ref{boundofnorms}
to estimate the norms at the places of $O_S$ and the standard
($p$--adic) Cauchy inequality for the the finite set  of  places
$S$.

\endssection
\endsection

\

\section Rational points on parabolic Riemann surfaces\par

\

Let $X_K$ be a projective variety of dimension $N>1$ defined over
some number field $K$ and $S\subset X_K(K)$. For every point $P\in
S$, we fix a  $LG$--germ $\V_P$ of dimension one on $P$. Suppose
that, there is a Riemann surface $M$ and an analytic map
$\gamma:M\to X_K(\Bbb C)$ (where we see $\Bbb C$ as a $K$--algebra
via the embedding $\sigma_0$ fixed in \S 1), passing through $S$.
Let $F\in M$ be a subset such that $\gamma (F)\subset S$. Suppose
that, for every $Q\in F$,  the germ of curve defined by $\gamma (M)$
near $\gamma(Q)$ coincides with $\V_{\gamma (Q)}$. In this chapter
we will suppose that $M$ is parabolic (cf.below) and the map
$\gamma$ is of finite order $\rho $ with respect to some positive
singularity $\tau$ on $M$. The aim of this chapter is to show that,
if the cardinality of $F$ is very big with respect to $\rho$, then
the image of $\gamma$ is algebraic. The theorem we will prove is
more general: it will concern points defined over extensions of
bounded degree.

\

\ssection Parabolic Riemann surfaces\par

\label parabolic. definition\par\defi A Riemann surface $M$ is said
to be {\it parabolic} if every upper bounded subharmonic function on
it is constant.\enddefi

\ssection Examples\par  a) The complex plane $\Bbb C$ is parabolic;

b) every algebraic Riemann surface is parabolic (affine or
projective);

c) if $F$ is a set of capacity zero on a parabolic Riemann surface
$M$ then $M\setminus F$ is again parabolic: so, for instance, the
complementary set of a lattice in $\Bbb C$ is parabolic;

d) a finite, ramified covering of a parabolic Riemann surface is
parabolic;

e) the unit disk {\it is not} parabolic.

f) Let $X$ be a {\it compact} Riemann surface and $f\colon Y\to X$
be a Galois covering with automorphism group $\Gamma$. It is proven
in [Gr] that $Y$ is parabolic if and only if $\Gamma$ has a subgroup
of finite index isomorphic to $\Bbb Z$ or $\Bbb Z\times \Bbb Z$ (or,
of course, if $\Gamma$ is finite).
\endssection
The notion of parabolic Riemann surface have been introduced, (and
studied), by Ahlfors (cf. [AS]); it is a class which strictly
contains the class of algebraic Riemann Surfaces; never the less,
Riemann surfaces in this class have many properties similar to the
algebraic. Unfortunately the name can be misleading: in the
uniformization theory of Riemann surfaces parabolic Riemann surfaces
are those having as universal covering the complex plane; here this
is not the case.
\smallskip
We recall here part of the theory and some of the properties of
parabolic Riemann surfaces; for the proofs see [AS] and [SN].

Let $M$ be a non compact parabolic Riemann surface and
$M^\ast=M\cup\{\infty\}$ be its Alexandroff compactification.

\label positivesingularity. definition\par\defi a) A {\rm positive
singularity} on $M$ is a couple $(\tau; U)$, where $U$ is a
neighborhood of $\infty$ (so, by definition, $M\setminus U$ is
compact) and $\tau$ is a positive harmonic function on $U$ such that

i) $\lim_{z\to\infty}\tau (z)=+\infty$;

ii) $\int_{\partial U}d^c\tau=-1$.

b) Two positive singularities $(\tau; U)$ and $(\tau'; U')$ are said
to be {\it equivalent} if $\tau -\tau'$ is a (upper and lower)
bounded harmonic function on $U\cap U'$.

We denote by $PS(M)$ the set of equivalence classes of positive
singularities on $M$.\enddefi

The relation between positive singularities and parabolic Riemann
surfaces is the following:

\prop The set $PS(M)$ is non empty.\endprop

It is even possible to prove that the fact that $PS(M)$ is non empty
is {\it equivalent} to the fact that $M$ is parabolic.

\ssection Examples\par a) If $M=\Bbb C$ then a positive singularity
on M is $\tau (z)=\log\vert z\vert^2$ (with $U=\{ \vert z\vert
>1\}$).

b) If $\pi: M\to\Bbb C$ is a finite ramified covering then
$\tau={{\pi^\ast(\log\vert z\vert^2)}\over{\deg\pi}}$ (and
$U=\pi^{-1}(\vert z\vert>1)$) is a positive singularity on $M$.
\endssection

Given a positive singularity $(\tau; U)$ on $M$ we can define the
Evans Kernel of it:

\label evanskernel. definition\par\defi  A $\tau$-- Evans kernel is
a function $e: M\times M\to (-\infty;\infty]$ such that

i) $-e(z; q)$ is subharmonic in $M$ as a function of $z$ and
harmonic in $M\setminus\{ q\}$;

ii) for every $q\in M$ we can find an open disk $\Delta_q$
neighborhood of $q$ such that
$$e(z;q)\vert_{\Delta_q}=-\log\vert z-q\vert^2 + \varphi$$
with $\varphi$ harmonic, $e(z;q)\geq 0$ and
$\int_{\partial\Delta_q}d^c_z(e(z;q))=-1$.

iii) for every fixed $q\in M$, the function $-e(z;q)$ is a positive
singularity on $M$ and it is equivalent to $\tau$;

iv) It is a symmetric function: $e(z;q)=e(q;z)$.
\enddefi

Remark that condition (iii) implies that, for a fixed $q$, we have
that $\lim_{z\to\infty}-e(z,q)=+\infty$. Moreover, for every couple
$q_1$ and $q_2\in M$ the two positive singularities $-e(z;q_i)$ are
equivalent; in particular there is a neighborhood of the infinity
where $\vert e(z;q_1)-e(z,q_2)\vert$ is uniformly bounded.

This definition will be useless without the:

\prop Given a positive singularity $\tau$ on $M$ there exists a
$\tau$--Evans Kernel unique up to a constant.
\endprop

For the proof see [SN] page 355 and ff.

We remark that the existence of the Evans kernel implies, in
particular,  that, given a positive singularity $\tau$ and a point
$q\in M$; there exists a positive singularity $\tau'$ equivalent to
$\tau$ and such that the implied open set is $M\setminus\{ q\}$.

Before we quote the regularity properties of the Evans kernel, we
recall the definition of the {\it Green functions}:

\label greenonsurfaces. definition\par\defi Let $U$ be a regular
region on a Riemann surface $M$ and $P\in U$. A {\rm Green function}
for $U$ and $P$ is a function $g_{U;P}(z)$ on $U$ such that:

a) $g_{U;P}(z)\vert_{\partial U}\equiv 0$ continuously;

b) $dd^cg_{U;P}=0$ on $U\setminus\{ P\}$;

c) near $P$, we have $g_{U;P}=-\log\vert z-P\vert^2+\varphi$, with
$\varphi$ continuous around $P$. \enddefi

We easily deduce from the definitions that
$dd^cg_{U;P}+\delta_P=\mu_{\partial U;P}$ where $\delta_P$ is the
Dirac at $P$ and  $\mu_{\partial U;P}$ is a positive measure of
total mass one and supported on $\partial U$.

Moreover the following is true:

\label uniquegreen. proposition\par\prop The Green function, if it
exists, it is unique.\endprop

We have the following regularity properties of the Evans kernels:

\label regularityofevans. proposition\par\prop a) The Evans Kernel
is a continuous map $e:M\times M\to (-\infty;\infty]$. Denote by
$\Delta\subset M\times M$ the diagonal; the Evans Kernel $e$ is
continuous and bounded on every open and relatively compact set
$U\subset M\times M$, such that $\overline U\cap\Delta=\emptyset$.

b) For every relatively compact open set $V\subset M$, we  have a
decomposition $e(p;q)=g_{V;q}(p)+v_{V;q}(p)$ with $v_{V;q}(\cdot)$
continuous and bounded.

c) Let $q_0\in M$ and $M_{\lambda}=\{z\in M\, /\,
e(z;q_0)\geq-\lambda\}$; then
$$e(p;q_0)=\lim_{\lambda\to\infty}(g_{M_{\lambda};q_0})-\lambda)$$

uniformly on every compact subset of $M\times M$, i.e. for every
$K\subset M$ compact, we have
$$\lim_{\lambda\to\infty}\sup_{(p;q)\in K^2}\vert
e(p;q)-(g_{M_{\lambda};q}-\lambda)\vert=0.$$
\endprop

\ssection Analytic maps of finite order from parabolic Riemann
surfaces\par

In this part we will describe the main definitions and properties of
the theory of analytic maps of finite order form a parabolic Riemann
surface to a projective variety. This theory is similar to the
classical Nevanlinna theory on maps from $\Bbb C$ to a projective
variety. Since we did not find an adequate reference we will give
some details.

We fix a parabolic Riemann surface $M$ and a positive singularity
$(\tau; U)$ on it (or more precisely its class in $PS(M)$). We also
fix a closed set $F\subset U$ such that:

\noindent i) the interior of $F$ is non empty;

\noindent ii) $\overline{M\setminus F}$ is compact.

Let $X$ be a projective variety and $\overline L$ be an ample line
bundle equipped with a $C^{\infty}$ positive metric and let $\gamma:
M\to X$ be an analytic map.

For every $t\in\Bbb R$ we define the {\it characteristic function}
$T_{\gamma,[\tau]}(t)$ of $\gamma$ with respect to the class
$[\tau]$ in the following way:
$$T_{\gamma; [\tau]}(t):=\int_F(t-\tau)^+\gamma^\ast(c_1(\overline
L)),$$

(where $(f(x))^+:=\sup\{ f(x);0\}$).

We will say that $\gamma$ has {\it finite order $\rho\in\Bbb R_{\geq
0}$} with respect to $[\tau]$ if
$$\limsup_{t\to\infty}{{\log T_{\gamma;[\tau]}(t)}\over{t}}=\rho.$$
This is equivalent to say that, for every $\epsilon>0$ and $t\gg 0$
we have $T_{\gamma; [\tau]}(t)\leq \exp(t(\rho+\epsilon))$. Making
the change of variables $t=\log r$, we see that $\gamma$ is of
finite order $\rho$ if, as soon as $r\to\infty$,  $T_{\gamma;
[\tau]}(\log r)\leq r^{\rho+\epsilon}$.

\ssection Example\par When $M=\Bbb C$ and $\tau(z)=\log\vert
z\vert^2$ then $T_{\gamma; [\tau]}(\log r)$ is the classical
characteristic function, defined, for instance, in [GK].\endssection

We prove now that the order of a map depends only on the map itself
and on the class of $\tau$ in $PS(M)$.  We introduce now the, so
called, non integrated version of the characteristic function:

Let $M$, $\tau$, $F$, $\gamma$ etc. as before. If  $t$ is
sufficiently big, we will introduce the following sets:  $F[t]:=\{
z\in F\,/\, \tau(z)\leq \log (t)\}$ and $B(r):=\{ (z;y)\in F\times
\Bbb R\,/\, \tau(z)\leq y\, {\rm and} \, y\leq\log (r)\}$. From the
definitions we obtain:
$$\eqalign{T_{\gamma; [\tau]}(\log(r))&=\int_F(\log
(r)-\tau(z))^+\gamma^\ast(c_1(\overline L))\cr
&=\int_{F[r]}(\int_{\tau(z)}^{log(r)}dy)\gamma^\ast(c_1(\overline
L)).\cr}$$

If we apply Fubini theorem we find that this last integral is
$$\eqalign{\int_{B(r)}dy\wedge p_1^\ast(\gamma^\ast(c_1(\overline
L)))&=\int_{-\infty}^{log
(r)}dy\int_{F[e^y]}\gamma^\ast(c_1(\overline L))\cr
&=\int_0^r{{dt}\over{t}}\int_{F[t]}\gamma^\ast(c_1(\overline
L)).\cr}$$ Call the function
$t_{\gamma;[\tau]}(y)=\int_{F[y]}\gamma^\ast(c_1(\overline L))$, the
{\it non integrated characteristic function} of $\gamma$ with
respect to $[\tau]$, and we obtain
$$T_{\gamma;
[\tau]}(\log(r))=\int_0^rt_{\gamma;[\tau]}(y){{dy}\over{y}}.$$

\ssection Example\par It is well known that,  if $M=\Bbb C$ and
$\tau =\log\vert z\vert^2$, then a possible definition of the
characteristic function is
$$T_{\gamma}(r)=\int_0^{r}{{dt}\over{t}}
\int_{\vert z\vert\leq t}\gamma^\ast(c_1(\overline L)).$$
\endssection

We are now ready to prove the independence of the order: \label
independenceoforder. theorem\par\thm The order $\rho$ of an analytic
map $\gamma$ from a parabolic Riemann surface $M$ to a projective
variety $X$ depends only on $\gamma$ and on the class $[\tau]\in
PS(M)$. more precisely it is independent on:

a) the choice of  the closed set $F$;

b) the choice of the representative $\tau\in [\tau]$;

c) the choice of the metrized ample line bundle $\overline L$ on
$X$.

\endthm

\Proof a) Independence on $F$:given two closed sets $F$ and $F'$,
then,  if the closure of the complementary of $F$ and $F'$ are
compact then also the closure of the complementary of $F\cap F'$ is
also compact; consequently we may suppose that $F'\subset F$. Let
$F'\subset F$ with $F$ and $F'$ as before. We denote by
$T^F_{\gamma}(t)$ ($T^{F'}_{\gamma}(t)$) and by $\rho^F$
($\rho^{F'}$ the characteristic function and the order respectively,
computed by using $F$ ($F'$). Since $T^F_{\gamma}(t)\geq
T^{F'}_{\gamma}(t)$, then $\rho^F\geq\rho^{F'}$.

On the other direction; denote
$\Vert\tau\Vert_{L^{\infty}(F\setminus F')}$ the $\sup$ of $\tau$
over $F\setminus F'$; remark that
$\Vert\tau\Vert_{L^{\infty}(F\setminus F')}<\infty$ because the
closure of $F\setminus F'$ is compact. Thus we have:
$$\eqalign{0\leq
T^{F}_{\gamma}(t)-T^{F'}_{\gamma}(t)&=\int_{F\setminus
F'}(t-\tau)^+\gamma^\ast(c_1(\overline L))\cr &\leq (\vert
t\vert+\Vert\tau\Vert_{L^{\infty}(F\setminus F')})\int_{F\setminus
F'}\gamma^\ast(c_1(\overline L))\cr}$$

Consequently we can find a constant $A\in\Bbb R_{\geq 0}$ such that
$T^{F}_{\gamma}(t)\leq T^{F'}_{\gamma}(t)+A\vert t\vert$. Suppose
$t$ (and consequently $T^{F}_{\gamma}(t)$ and $T^{F'}_{\gamma}(t)$)
sufficiently big (otherwise $\rho^F$ and $\rho^{F'}$ are both zero);
we have
$$\eqalign{\rho^F=\limsup_{t\to\infty}{{\log(T^{F}_{\gamma}(t))}\over{t}}
&\leq\limsup_{t\to\infty}{{\log(T^{F'}_{\gamma}(t)+A\vert
t\vert)}\over{t}}\cr
&\leq\limsup_{t\to\infty}{{\log(T^{F'}_{\gamma}(t))+\log(A\vert
t\vert)}\over{t}}=\rho^{F'}.\cr}$$

b) Independence on $\tau\in [\tau]$: First of all it is easy to see
that if $\tau'=\tau +A$ (for a suitable constant $A\in \Bbb R$),
then, if we compute the order of $\gamma$  by using $\tau'$ or by
using $\tau$ we obtain the same number. In general, if $\tau_1$ and
$\tau_2$ are in the same class in $PS(M)$, then, over suitable open
sets, we can suppose that there exist two constants $A$ and $B$ such
that $A\leq\tau_2-\tau_1\leq B$. Consequently we can suppose that
$\tau_1\leq\tau_2$. Put $F_i[y]:=\{ \tau_i(z)\leq\log y\}$
($i=1,2$). Then,  for $y\gg 0$, $F_2[y]\subseteq F_1[y]$,
consequently, one easily see, by using the non integrated form of
the characteristic function, that $T_{\gamma,\tau_1}(t)\geq
T_{\gamma,\tau_2}(t)$. From this we conclude.

c) Independence on $\overline L$: First of all we check that the
definition of the order is independent on the choice of the metric
on the ample line bundle $L$: Suppose that we choose two (positive)
metrics $\Vert\cdot\Vert_1$ and
$\Vert\cdot\Vert_2=\Vert\cdot\Vert_1\cdot\exp(\varphi)$ on $L$, with
$\varphi$ a bounded $C^{\infty}$ function on $X$. If we denote by
$T_{\gamma,i}(t)$ the characteristic function computed by using the
metric $\Vert\cdot\Vert_i$, we have
$$T_{\gamma,2}(\log(r))=T_{\gamma,1}(\log(r))+
\int_0^r{{dy}\over{y}}\int_{F[y]}\gamma^\ast dd^c\varphi.$$ By
Stokes Theorem we have
$$\eqalign{\int_0^r{{dy}\over{y}}\int_{F[y]}\gamma^\ast dd^c\varphi&=
\int_0^r{{dy}\over{y}}\int_{\partial F[y]}\gamma^\ast d^c\varphi\cr
&=\int_{F[e^r]}d\tau\wedge\gamma^\ast d^c\varphi.\cr}$$ Thus since
$\tau$ is harmonic, the following equality holds:
$d\tau\wedge\gamma^\ast d^c\varphi=d\gamma^\ast\wedge
d^c\tau=d(\gamma^\ast\varphi \cdot d^c\tau)$. By applying Stokes
again, we obtain that the last integral is equal to
$$\int_{\partial F[e^r]}\gamma^\ast\varphi\cdot d^c\tau.$$
We deduce the independence on the metrics because of property (ii)
of \ref{positivesingularity} and the fact that $\varphi$ is bounded.

In order to prove that the order of $\gamma$ is independent on the
choice of the ample line bundle $L$, it is sufficient to remark
that, if $L_1$ and $L_2$ are two ample line bundles on $X$, such
that $L_1\otimes L_2^{-1}$ is ample then, for $t\gg 0$,
$T_{\gamma;L_1}(t)\geq T_{\gamma; L_2}(t)$, and, for every positive
integer $D$ we have $T_{\gamma; L_1^D}(t)=DT_{\gamma; L_1}(t)$.

\

\ssection Points on the image of maps of finite order\par

Let $M$ be a (non compact) parabolic Riemann surface and $\tau$ a
positive singularity on $M$. Let $e(p;q)$ be the $\tau$--Evans
Kernel over $M$. We can  define a "canonical" metric on the tangent
bundle $TM$ and, for every point $p\in M$ on the line bundle
$\O_M(p)$: Let $\Delta$ be the diagonal divisor on $M\times M$; we
put a metric on the line bundle $\O_{M\times M}(\Delta )$ in the
following way: if ${\Bbb I}_{\Delta}$ is the section of $\O_{M\times
M}(\Delta )$ defining $\Delta$, we define $\Vert{\Bbb
I}_{\Delta}\Vert(p;q):=\exp(-e(p;q))$. Observe that, by property
(ii) of the Evans Kernel (definition \ref{evanskernel}), the metric
on $\O_{M\times M}(\Delta )$ is well defined.

By pull--back, we define a metric on $TM:=\Delta^\ast(\O_{M\times
M}(\Delta))$; for every point $p\in M$ let $\beta_p:M\to M\times M$
is the map $\beta_p(q)=(p;q)$;  on $\O_M(p)$ we put the metric
induced by the isomorphism $\O_M(p)=\beta_p^\ast(\O_{M\times
M}(\Delta))$.

\rmk If we choose the "other'' embedding $\alpha_p(q):=(q;p)$ we
obtain the same metric on $\O_M(p)$ because of the symmetry of the
Evans Kernel.
\endrmk

By construction, for every point $p\in M$ we find a canonical
adjunction isometry
$$TM\vert_p\simeq\O_M(p)\vert_p.$$

Let $X$ be a projective variety, $\overline L$  a metrized ample
line bundle over $X$ (with positive metric) and $\gamma:M\to X$ be
an analytic map of finite order $\rho$ with respect to $[\tau]$. We
suppose that $\gamma (M)$ is Zariski dense in X. For every positive
integer, we have an injective map
$$\gamma^\ast:H^0(X;L^D)\hookrightarrow H^0(M;\gamma^\ast(L^D)).$$
We fix a finite set of points $F:=\{ p_1;\dots ;p_s\}\subset M$.
and, for every positive integer $i$, we consider the canonical
injective map
$$\eta_i:H^0(M;\gamma^\ast L^D(-i\sum_{p_j\in
F}p_j))\hookrightarrow H^0(M;\gamma^\ast L^D).$$

We denote by $E_D^i$ the subspace
$(\gamma^\ast)^{-1}(\eta_i(H^0(M;\gamma^\ast L^D(-i\sum_{p_j\in
F}p_j))))$ and by $\gamma_D^i$ the composite map
$$\eqalign{&E_D^i\buildrel{\gamma^\ast}\over\longrightarrow
\eta_i(H^0(M;\gamma^\ast L^D(-i\sum_{p_j\in
F}p_j)))\cr&\longrightarrow \eta_i(H^0(M;\gamma^\ast
L^D(-i\sum_{p_j\in F}p_j)))\big/ \eta_{i+1}(H^0(M;\gamma^\ast
L^D(-(i+1)\sum_{p_j\in F}p_j)))\cr &\simeq\bigoplus_{p_j\in
F}(\gamma^\ast(L^D)\otimes (TM)^{-i})\vert_{p_j};\cr}$$ the
$E_D^i$'s are finite dimensional vector spaces equipped with the
$\sup$ norm. The vector space $\bigoplus_{p_j\in
F}(\gamma^\ast(L^D)\otimes (TM)^{-i})\vert_{p_j}$ is also a finite
dimensional hermitian vector space: we equip it with the "direct
sum'' metric. Eventually, if $p_h\in F$, we let
$\gamma_{D,h}^i:E^i_D\to\gamma^\ast(L^D)\otimes
(TM)^{-i})\vert_{p_h}$  be the linear map composite of $\gamma_D^i$
with the canonical projection $\bigoplus_{p_j\in
F}(\gamma^\ast(L^D)\otimes
(TM)^{-i})\vert_{p_j}\to\gamma^\ast(L^D)\otimes
(TM)^{-i})\vert_{p_h}$.

The main theorem of this chapter is the following: if, in Theorem
\ref{generalalgebraization}, the involved analytic germs come from a
parabolic Riemann surface and the corresponding map is of finite
order, then we can prove strong version of Schwartz lemma. More
precisely:

\label schwartzforparabolic. theorem\par\thm Let $\epsilon>0$. Under
the hypotheses above, there is a constant $C$ depending only on
$\gamma$, on $F$ etc., but independent on $i$ and $D$ (provided that
${{i}\over{D}}$ is sufficiently big), such that, for every $h$, we
have
$$\log\Vert\gamma_{D,h}^i\Vert\leq-{{i\cdot
Card(F)}\over{\rho+\epsilon}}\log({{i}\over{D}})+i\cdot C.$$
\endthm

\Proof Let $p:=p_h$. Let $r$ be a sufficiently big positive real
number such that
$$\Omega_r:=\left\{ z\in M\, /\, \Vert {\Bbb I}_p\Vert \leq
r\right\}\supset F;$$ thus, by definition $\Omega_r=\{ z\in M\,/\,
-e(z;p)<\log (r)\}$. On $\Omega_r$ we consider the function
$g_{r;p}(z):=\log (r)+e(z;p)$; this function has the following
properties:

-- near $p$ it is of the form $-\log\vert z-p\vert^2+\varphi_p(z)$
with $\varphi_p(z)$ harmonic near $p$;

-- $\lim_{z\to\partial\Omega_r}g_{r;p}=0$ uniformly;

-- $g_{r;p}$ is continuous and harmonic in $\Omega_r\setminus\{
p\}$;

-- it is positive.

If we define $g_{r;p}\equiv 0$ on $M\setminus\Omega_r$  we can
extend $g_{r;p}$ to a continuous function on $M$, which we will
denote again by $g_{r;p}$. By construction $g_{r;p}$ is the Green
function of $p$ in $\Omega_r$ (cf. \ref{uniquegreen}). In particular
$dd^cg_{r;p}=-\delta_p+\mu_{p;\partial\Omega_r}$ (as distributions),
where $\mu_{p;\partial\Omega_r}$ is a positive measure of total mass
one and supported on $\partial\Omega_r$.

Let $s\in E_D^i$; thus, by definition,
$\gamma^\ast(s)=\eta_i(\tilde{s})$, where $\tilde s$  is  a global
section of $\gamma^\ast(L^D)(-i\sum P_j)$. We will denote by $\Vert
s\Vert$ the norm of $s$ as section of the hermitian line bundle
$\gamma^\ast(L^D)$ and by $\Vert\tilde s\Vert$ the norm of $\tilde
s$ as section of the hermitian line bundle $\gamma^\ast(L^D)(-i\sum
P_j)$. By  Stokes Theorem we have
$$\int_M\log\Vert\tilde s\Vert^2\cdot
dd^cg_{r;p}=\int_Mdd^c\log\Vert\tilde s\Vert^2\cdot g_{r;p}.$$
Consequently we get
$$\int_M\log\Vert\tilde
s\Vert^2(-\delta_p+\mu_{p;\partial\Omega_r})=\int_M(\delta_{div(\tilde
s)}-D\gamma^\ast (c_1(\overline L))\cdot g_{r;p};$$ where
$\delta_{div (\tilde s)}$ is the Dirac measure supported on the
effective "divisor" of the zeros of $\tilde s$. Since
$\mu_{p;\partial\Omega_r}$ is positive of total mass one and
$g_{r;p}$ is positive we obtain
$$\eqalign{-\log\Vert\tilde s\Vert^2(p)+\log\Vert
s\Vert^2_{L^{\infty}}+&i\sum_{p_j\in
F}\int_Me(\cdot;p_j)\mu_{p;\partial\Omega_r}\cr &\geq
-D\int_M\gamma^\ast(c_1(L))\cdot (\log r+e(\cdot ;p))^+.\cr}$$ By
property (iii) of the Evans Kernel, (cf. \ref{evanskernel}) we then
obtain
$$\log\Vert s\Vert^2_{L^{\infty}}+i\sum_{p_j\in F}\int_Me(\cdot
;p_j)\mu_{p;\partial\Omega_r} + DT_{\gamma;[\tau]}(\log
(r))\geq\log\Vert \tilde s\Vert^2(p).$$

By applying Prop \ref{regularityofevans} there exists a constant $B$
and an open set $U\subseteq M$ such that $M\setminus U$ is compact,
contains $F$ and, for every $p_j\in F$ and for $z\in U$ we have
$\vert e(z;p_j)-e(z;p)\vert\leq B$ (they both define equivalent
positive singularities). So on $\partial\Omega_r$ (which, for $r\gg
0$ we can suppose contained in $U$) we have
$\vert\log(r)+e(z;p_j)\vert\leq B$; from this we get
$$\log\Vert s\Vert^2_{L^{\infty}}-i\cdot Card (F)\log (r)+ DT_{\gamma;[\tau]}(\log
(r)) +i\cdot B_1\geq\log\Vert \tilde s\Vert^2(p).$$

By hypothesis $\gamma$ is of finite order $\rho$, thus, for $r\gg
0$, we have $T_{\gamma ;[\tau]}(\log (r))\leq r^{\rho + \epsilon}$,
consequently
$$\log\Vert s\Vert_{L^{\infty}}-i\cdot Card (F)\log (r)+ Dr^{\rho
+\epsilon}+i\cdot B_1\geq\log\Vert \tilde s\Vert(p).$$

Let $h(r)=-i\cdot Card (F)\log (r)+ Dr^{\rho +\epsilon}$; we have
${{dh}\over{dr}}(r)=-{{i\cdot Card
(F)}\over{r}}+(\rho+\epsilon)Dr^{\rho+\epsilon -1}$; thus $h(r)$ has
a minimum in $r_0=\left({ {i\cdot Card
(F)}\over{(\rho+\epsilon)D}}\right)^{{{1}\over{\rho+\epsilon}}}$,
which, as soon as ${{i}\over{D}}$ is sufficiently big is allowed. In
this case
$$h(r_0)=-{{i\cdot Card (F)}\over{\rho+\epsilon}}\log({{i}\over{D}})+
i\cdot B_2$$ with $B_2$ independent on $i$ and $D$. From this we
conclude.

\smallskip

Now we can state, and prove the main theorem of this section. As in
the introduction, we suppose that $K$ is a number field and
$\sigma_0:K\hookrightarrow\Bbb C$ is an embedding of $K$ in $\Bbb
C$. We also fix an embedding of the algebraic closure $\overline K$
of $K$ in $\Bbb C$. Let $X$ be a quasi projective variety of
dimension $N$ defined over $K$; Let $S\subseteq X(\overline K)$ and,
for every positive integer $r$ denote by $S_r$ the set $\{ x\in S\;
s. t.\; [\Bbb Q(x):K]\leq r\}$.

\label maintheoremforriemannsurfaces. theorem\par\thm Let $M$ be a
parabolic Riemann surface (with a fixed positive singularity). Let
$\gamma\colon M\to X(\Bbb C)$ be an holomorphic map of finite order
$\rho$ with Zariski dense image. Suppose that, for every $\overline
K$--rational point $P\in S\cap\gamma (M)$, the formal germ $\hat
M_P$, of $M$ near $P$, is (the pull back of) a $LG$-- germ of type
$\alpha$ (in its field of definition). then
$${{Card (\gamma^{-1}(S_r))}\over{r}}\leq {{N+1}\over{N-1}}\rho\alpha[K:\Bbb Q].$$
\endthm

\Proof It suffices to prove that, given a finite set of points $F=\{
P_1,\dots, P_s\}\subset \gamma^{-1}(S_r)$ and $\epsilon>0$, we have
$${{Card(F)}\over{r}}\leq {{N+1}\over{N-1}}(\rho+\epsilon)\alpha[K:\Bbb Q].$$
Let $F$ be such a set and $L$ be a finite extension of $K$ which
contains $\Bbb Q(\gamma(P_j))$ for every $P_j\in F$.

For every number field $L$ we will denote by $M_L$ the set of places
of $L$. We  fix $P_j$ and, and let $\sigma_{P_j}\in S_{\Bbb Q(P_j)}$
be the embedding extending $\sigma_0$. If $i/D\gg 0$, from Corollary
\ref{finitenorms}  and Theorem \ref{schwartzforparabolic} we get,

$$\eqalign{{{1}\over{[\Bbb Q(P_j):\Bbb Q]}}&\sum_{\sigma\in M_{\Bbb
Q(P_j)}}\log\Vert \gamma_{D,j}^i\Vert_\sigma\cr &={{1}\over{[\Bbb
Q(P_j):\Bbb Q]}}\left( \log\Vert
\gamma_{D,j}^i\Vert_{\sigma_{\sigma_{P_j}}}+\sum_{\sigma\neq\sigma_{P_j}}\log\Vert
\gamma_{D,j}^i\Vert_\sigma\right)\cr &\leq
-{{Card(F)}\over{(\rho+\epsilon)\cdot[\Bbb Q(P_j):\Bbb
Q]}}i\log\left({{i}\over{D}}\right)+\alpha\cdot i\log
(i)+C_2\left(i+D)\right),\cr}$$ with $C_2$ independent on $i$, $D$
and $\Bbb Q(P_j)$. Thus, we deduce
$$\eqalign{{{1}\over{[L:\Bbb Q]}}&\sum_{\sigma\in
M_L}\log\Vert\gamma_D^i\Vert_\sigma\leq{{1}\over{[L:\Bbb
Q]}}\sum_{\sigma\in M_L}\sup_{P_j\in F}\log\Vert
\gamma_{D,j}^i\Vert_\sigma+\log(Card(F))\cr &\leq
-{{Card(F)}\over{r\cdot[K:\Bbb
Q](\rho+\epsilon)}}i\log\left({{i}\over{D}}\right)+\alpha\cdot i\log
(i)+C_2\left(i+D\right)+\log(Card(F)).\cr}$$ But since the image is
Zariski dense and the dimension of $X$ is $N$, we can apply Theorem
\ref{generalalgebraization} and we obtain that
$${{Card(F)}\over{r\cdot[K:\Bbb Q]}}\leq {{N+1}\over{N-1}}(\rho+\epsilon)\alpha.$$
and we eventually conclude.

\smallskip

We state also the following corollaries, which show the the link
with the classical Schneider--Lang Theorem (where $X=\Bbb C^n$ and
$M=\Bbb C$):

\label slonqbar. corollary\par\cor Suppose we are in the hypotheses
above, then
$$\sum_{P\in\gamma^{-1}(S)}{{1}\over{[\Bbb Q(P);
K]}}\leq {{N+1}\over{N-1}}\rho\alpha.$$
\endcor

If $X$ is a quasi projective variety defined over the number field
$K$ and $r$ is a positive integer, denote by $X_r$ the set $\{ P\in
X(\overline K)\;\; {\rm s.t.}\;\; [\Bbb Q(P):K]\leq r\}$.

\label classicalsl. corollary\par\cor Let $X$ be an algebraic
variety defined over a number field $K$ and let $F\hookrightarrow
T_X$ be a  foliation of rank one (defined over K). Suppose that the
holomorphic foliation $F_{\sigma}\subset (T_X)_{\sigma}$ has a leaf
$M$ which is parabolic of finite order $\rho$ (for some positive
singularity on $M$) whose Zariski closure has dimension $d>1$, then
$${{Card((X_r\setminus Sing(F))\cap M)}\over{r}}\leq {{d+1}\over{d-1}}\rho[K:\Bbb Q].$$
\endcor

\label bostgk. corollary\par\cor If, in the hypotheses of
\ref{classicalsl}, the foliation is closed under $p$ derivation for
almost all primes $\p$ of $O_K$, then the leaf  passing through an algebraic point is an algebraic
curve.\endcor

The proofs of the corollaries are straightforward applications of
the main theorem.

\rmk One should notice that corollary \ref{bostgk} is a particular
case of the main theorem of [Bo].\endrmk

\endssection
\endsection

\

\

\section Rational points on higher dimensional subvarieties\par

\

In this section we will deal with Kh\"aler varieties which enjoy
properties similar to the parabolic Riemann surfaces. We will study
analytic maps between these varieties and Projective varieties. When
one look carefully to the previous section one sees that the main
tool to prove the main theorems was the existence of the Evans
kernel on the Riemann surface we are working with. We will see that
The analogous of a Evans kernel will suffice to develop a value
distribution theory and the notion of the order of an analytic map.

\label laplaceparab. definition\par\defi Let $A$ be a $d$
dimensional Kh\"aler manifold with a fixed Kh\"aler metric $\omega$.
We will say that $(A,\omega)$ (or simply $A$) is {\rm conformally
parabolic} if there exists a function
$$g:A\times A\longrightarrow (-\infty;\infty]$$
with the following properties:

a) For every $p\in A$ the function $g_p(z):=g(p,z)$ is $C^{\infty}$
in $A\setminus\{ p\}$ and satisfy the following differential
equation:
$$dd^c(g_p)\wedge\omega^{d-1}=\delta_p;$$
where $\delta_p$ is the Dirac measure concentrated in  $p$;

b) for every $p$ we have that
$$\lim_{z\to\infty}g_p(z)=+\infty;$$

c) for every couple $p$ and $q$ there exists a neighborhood of the
infinity $U$ and a constant $C$ such that, for every $z\in U$
$$\vert g_p(z)-g_q(z)\vert\leq C.$$

\enddefi

One easily sees that if $d=1$, up to the sign, $g$ is an Evans
kernel.

We will call the function $g$, a {\it Evans Kernel} of the
conformally parabolic manifold $A$.

Before we develop a value distribution theory on conformally
parabolic varieties, we give the main example we have in mind.

\ssection Examples of conformally Parabolic varieties\par Let
$\overline A$ be a compact Kh\"aler manifold of dimension $d$,
$\omega$ a Kh\"aler form on it and $H_\infty$ be a divisor on
$\overline A$. Let  $A$ be the open set $\overline A\setminus
H_\infty$. We will denote by $R$ the degree of $H_\infty$ with
respect to $\omega$, we also fix a $C^\infty$ metric
$\Vert\cdot\Vert_{H_\infty}$ on $\O_{\overline A}(H_\infty)$.

Let $p\in A$ be a point.

\label functiongp. theorem\par\thm There exists a function
$g_p\colon \overline A\to[-\infty,+\infty]$ unique up to an additive
constant with the following properties:

\item{(i)} $g_p$ is $C^{\infty}$ outside $\{p\}\cup H_{\infty}$;

\item{(ii)} $g_p$ is a solution of the following differential
equation:
$$dd^c(g_p)\wedge\omega^{d-1}=\
\delta_p-{{1}\over{R}}\delta_{H_{\infty}}\omega^{d-1};$$ $\delta_p$
(resp. $\delta_{H_{\infty}}$) being the Dirac measure concentrated
on $p$ (resp. $H_{\infty}$).
\endthm

We will only sketch the proof of \ref{functiongp}. Details can be
filled by standard Hodge theory, cf. for instance [Vs].

\Proof Let $T$ be the current
$\delta_p-{{1}\over{R}}\delta_{H_{\infty}}\omega^{d-1}$. It is a
current of bi--degree $(d,d)$ and $\int_{\overline A}T=0$.
Consequently, by the Hodge decomposition, there exists a current $S$
such that $\Delta_{\overline{\partial}}(S)=T$
($\Delta_{\overline{\partial}}$ being the laplacian with respect to
the metric $\omega$). Since $T$ is $C^{\infty}$ outside $\{P\}\cup
H_{\infty}$, one can prove, by following the proof of 6.32 chapter 6
of [Wa] that also the current $S$ is $C^{\infty}$ outside $\{P\}\cup
H_{\infty}$. Let $L$ be the operator on forms obtained by wedging
with the form $\omega$. Denoting by ${\cal D}^{i,i}$ the sheaf of
currents of bidegree $(i,i)$, it is well known that $L^d$ induces an
isomorphism between ${\cal D}^{0,0}$ and ${\cal D}^{d,d}$. Thus,
there is a distribution $\tilde g_p$ such that $L^d(\tilde g_p)=S$.
By the standard commutation rules between $L$,
$\Delta_{\overline{\partial}}$, $\partial$, $\overline{\partial}$
and $\overline{\partial}^{\ast}$, we can find a constant $c$ (which
depends only on $d$) such that, for every local section $g\in{\cal
A}^{0,0}$,
$$dd^c(g)\wedge\omega^{d-1}=c\cdot\Delta_{\overline\partial}(g)\wedge\omega^d=c\cdot
L^d(\Delta_{\overline\partial}(g)).$$ Moreover
$L^d(\Delta_{\overline\partial}(g))=\Delta_{\overline\partial}(L^d(g))$.
Consequently we can find a function $g_p$ with the properties
claimed by the statement. The difference of two functions verifying
(i) and (ii) will be an harmonic function on $\overline A$, thus it
will be a constant.

\label evansofkhaler. theorem\par\thm The function $g(p,z):=g_p(z)$
is an Evans kernel on $A$.
\endthm

\Proof Let $s_{\infty}\in H^0(\UA,\O_{\UA}(H_\infty))$ be a section
such that $div(s_{\infty})=H_{\infty}$; then we claim that there
exists a $C^{\infty}$  function $f_p$ on $\UA\setminus\{p\}$ such
that, on $\UA\setminus\{p\}$,
$$g_p=-{{1}\over{R}}\log\Vert s_{\infty}\Vert^2_{H_\infty}+f_p.$$
Indeed by Poincar\'e--Lelong formula,  $dd^c\log\Vert
s_{\infty}\Vert^2_{H_\infty}=\delta_{H_{\infty}}-c_1(\O_{\UA}(H_\infty);\Vert\cdot\Vert_{H_\infty})$,
we have that, on $\UA\setminus\{p\}$
$$dd^c(g_p+{{1}\over{R}}\log\Vert
s_{\infty}\Vert^2_{H_\infty})\wedge\omega^{d-1}=-{{1}\over{R}}c_1(\O_{\UA}(H_\infty);\Vert\cdot\Vert_{H_\infty})\wedge\omega^{d-1};$$
part the claim follows by the same method of the proof of
\ref{functiongp}. Consequently properties (a), (b) and (c) are
easily verified.

\label khalerlaplaceparabolic. corollary\par\cor The complementary
of an effective divisor of a compact Kh\"aler manifold is
conformally parabolic. In particular every quasiprojective variety
is conformally pa\-ra\-bo\-lic.
\endcor

\endssection

\ssection counting functions\par

In this subsection we will develop a value distribution theory for
conformally parabolic manifold which is analogue to the value
distribution theory on parabolic Riemann surfaces.

We fix a $d$ dimensional ($d\geq 2)$) conformally parabolic manifold
$A$ with Kh\"aler form $\omega$ and Evans kernel $g$ (we suppose
that $d\geq 2$ because the case $n=1$ is treated in the previous
section and it is a little bit different on the estimates).

Before we need to  state and prove some of the properties of the
functions $g_p$'s.

\label propertiesofgp. proposition\par\prop Let $p\in A$. There
exists a neighborhood $U_p$ of $p$, analytically equivalent to a
$d$--th dimensional ball centered at $p$ (with coordinates $z$) and
a continuous non vanishing function $\chi_p$, such that
$$\vert g_p\vert\vert_{U_p}\cdot\Vert z\Vert^{2(d-1)}=\chi_p(z).$$
Moreover there is a positive constant $\alpha$ such that
$$\chi_p(z)={{1}\over{n-1}}+O_{\omega}(\Vert z\Vert^\alpha).
$$

\endprop
\Proof the stated properties can be deduced from Theorem 4.13 page
108 of [Au]; more precisely from formula (17) page 109, Lemma 4.12
page 107 and formula (8) page 106 of loc. cit.

\rmk One can deduce from the proof of Theorem 4.13 page 109 of [Au]
that $\chi_p(z)$ is a function which is regular enough for our
purposes. Indeed it will be of class $C^{1,\alpha}$ for every
$\alpha<1$.
\endrmk

For every continuous function $\lambda: A\to [-\infty, +\infty)$ we
define the following sets:

$$S_\lambda(t):=\left\{ z\in A \; {\rm such \; that }\;
\lambda(z)=t\right\};$$ and $$B_\lambda(t):=\left\{ z\in A \; {\rm
such \; that }\; \lambda(z)\leq t\right\}.$$

Moreover, if $D$ is a subset of $A$, we denote by $D_{\lambda}(t)$
the set $D\cap B_{\lambda}(t)$.

Remark that, by condition (b)of the definition of the Evans kernel,
if $p\in A$, $S_{g_p}$ and $B_{g_p}$ are compact.

If $\alpha$ is a $(1,1)$ form on $A$ we define the {\it counting
function} with respect to $\alpha$ to be
$$\eqalign{T_{\alpha}(r)&:=\int_0^r{{dt}\over{t}}\int_{B_{g_p}(\log(t))}\alpha\wedge(\omega)^{d-1}\cr
&=\int_{-\infty}^{\log(r)}dt\int_{g_p\leq
t}\alpha\wedge(\omega)^{d-1}\cr &
=\int_A(\log(r)-g_p)^+(\alpha\wedge (\omega)^{d-1})\cr}$$ (where, as
in \S 4, $(f)^+$ means $\sup\{ f, 0\}$).

Let $X$ be projective variety and $\gamma\colon A\to X$ be an
analytic map. Suppose that $\overline L$ is an ample line bundle
equipped with a positive metric. We define $T_{\gamma}(r)$ to be
$T_{\gamma^{\ast}(c_1(\overline L))}(r)$.

\label finiteorderofaffine. definition\par\defi We say that
$\gamma\colon A\to X$ is {\it of finite order} $\rho$ if
$$\limsup_{r\to\infty}{{\log(T_\gamma(r))}\over{\log(r)}}=\rho.$$
\enddefi

As in the previous section, this means that, for every $\epsilon
>0$ there is a $r_0$ such that for $r\geq r_0$ we have
$$T_\gamma(r)\leq r^{\rho +\epsilon}.$$

The order depends only on the map $\gamma$; in particular it is
independent on the choice of the ample line bundle $L$ (and on the
metric on it).

\label orderiswelldefined. proposition\par\prop The order of
$\gamma: A\to X$ depends only on $\gamma$; more precisely:
\item{a)} the order is independent on the choice of the point $p$;
\item{b)} the order is independent on the choice of the ample line
bundle with positive metric $L$  on $X$.
\endprop
\Proof Let $p_1$ and $p_2$ be two points on $A$. By  property (c) of
the Evans kernel, we can find a compact set $K$ containing  the
$p_i$'s and constants $C_i$ such that, outside $K$, $g_{p_2}+A_2\leq
g_{p_1}\leq g_{p_2}+C_1$. Consequently the proof is analogous to the
case of parabolic Riemann surfaces and we leave the details to the
reader.

We will show now that we can use the counting function $T(r)$ in
order to control the norm of jets of sections of line bundles.

Let $L$ be an hermitian line bundle equipped with a positive metric
on $A$. Let $s\in H^0(A,L)$. We define the {\it proximity function}
of $s$ to be
$$\mu_p(\log\Vert s\Vert^2)(r):=\int_{S_{g_p}(\log(r))}\log\Vert s\Vert^2\cdot
d^c(g_p)\wedge\omega^{d-1}.$$ By Stokes theorem we have
$$\int_{S_{g_p}(\log(r))}d^c(g_p)\wedge\omega^{d-1}=1$$ and, for
$r\gg 0$, $d^cg_p\wedge\omega^{d-1}$ is a positive measure on
$S_{g_p}(r)$. Consequently, if $\Vert s\Vert_{\infty}<\infty$ and
$r\gg 0$, \labelf boundonlinfinity\par$$\log\Vert
s\Vert^2_{\infty}\geq\mu_p(\log\Vert
s\Vert^2)(r).\eqno{{(\numfo)}}$$\advance\ssnu by1

The following theorem is an analogue of the Nevanlinna F.M.T., in
this contest.

\label fmt. proposition\par\prop For every $R$ and $r$ with $R>r$
the following equality holds:
$$\mu_p(\log\Vert s\Vert^2)(R)-\mu_p(\log\Vert
s\Vert^2)(r)=\int_{\log(r)}^{\log(R)}dt\int_{g_p\leq
t}(\delta_{div(s)}-c_1(L))\wedge(\omega)^{d-1}.$$
\endprop

The proof is a direct application of Stokes formula as in the proof
of the Nevanlinna F.M.T. (cf [GK]).

Let $L$ be a line bundle on $A$ equipped with a positive metric. Let
$s\in H^0(A,L^D)$ be a global section of it. Suppose that $s$
vanishes at the order $i$ in $p$. Let $j^i(s)$ be its $i$--th jet in
$p$. It is a well defined section of $S^i(T_A^{\ast})\otimes
L^D\vert_p$. Since we fixed a metric on $A$ and on $L$, the vector
space $S^i(T_A^{\ast})\otimes L^D\vert_p$ is equipped with an
hermitian metric; we will denote by $\Vert\cdot\Vert_{i,D}$ the
induced norm on it.

We will give a bound of the norm of $j^i(s)$ in terms of
$\mu_p(\log\Vert s\Vert^2)(R)$. This can be seen as an analogue of
the classical Cauchy inequality on the complex plane.

\label generalboundofnorms. proposition\par\prop There exists a
constant $C$ depending only on $r$ and on the metric on $L$ such
that, for every $s\in H^0(A,L^D)$ as above, we have
$$\mu_p(\log\Vert s\Vert^2)(r)\geq \log\Vert
j^i(s)\Vert^2_{i,D}+C(i+D).$$
\endprop

Since this statement do not change if we change the metric on the
tangent bundle $T_A$ (but if we change it of course, the constant
$C$ will change), we are free to choose the metric we prefer on the
vector space $T_A\vert_p$.

We consider the blow up $\tilde A$ of $A$ in $p$ let $E$ be the
exceptional divisor. We can define a metric on ${\cal O}(E)$ in the
following way: $\Vert E\Vert (z):={{1}\over{\vert
g_p(z)\vert^{1/2(d-1)}}}$. The Kahler metric $\omega$ induces a
Fubini--Study  metric $\omega_{FS}$ on $E$. Let $\tilde s$ be the
strict transform of $s$; it is a section of the hermitian line
bundle $L(-iE)$. We know  that $\log\Vert
 j^i(s)\Vert^2_{i,D}=\int_E\log\Vert\tilde s\Vert^2\omega_{FS}^{d-1}+C$
where $C$ is a controlled constant (cf. [Bo] \S 4.3.2). By Theorem
4.13 page 109 of [Au] one can easily show that
$$\lim_{R\to-\infty}\int_{g_p=R}\log\Vert\tilde s\Vert^2
d^cg_p\wedge\omega^{d-1}=\int_E\log\Vert\tilde
s\Vert^2\omega_{FS}^{d-1}.$$

Proposition \ref{generalboundofnorms} will be consequence of the
following more general statement.

\label generalboundofnonrms2. theorem\par\thm There exists a $\tilde
r\gg 0$ and a $\gamma>0$ depending only on the Kahler metric
$\omega$ such that, if $\log(r)>\tilde r$ the following holds:
$$\eqalign{ T_{c_1(L)}(r)&\geq -\mu_p(\log\Vert s\Vert^2)(r)+\log\Vert
j^i(s)\Vert_{i,D}\cr &+\int_{-\tilde
r}^{\log(r)}dt\int_{div(s)_{g_p}(t)}\omega^{d-1}-{{i}\over{d-1}}\cdot\left(
\log(\tilde r)+O_{\omega}({{1}\over{\tilde
r^{\gamma}}})\right).\cr}$$
\endthm

\Proof We start with the following: Since the metric $\omega$ is
K\"ahler, we can choose a neighborhood $U_p$ of $p$, analytically
equivalent to the unit ball and coordinates $(z_1,\dots z_d)$
centered at $p$ on it, in such a way that the following property
holds: denote by $\lambda=dd^c\vert z\vert^2$ the standard Euclidean
(1,1) form on $U_p$; then $\omega\vert_{B_r}=\lambda+k$ with
$k=O(\vert z\vert^2)$ when $\vert z\vert\to 0$.

\label euclideanfact. lemma\par\lemma Suppose that $0<\tilde r<R$
and $\tilde r$ is such that $B_{g_p}(-\tilde r)\subseteq U_p$; then
we can find a positive constant $\gamma$ such that
$$\int_{-R}^{-\tilde r}dt\int_{div(s)_{g_p}(t)}\omega^{d-1}\geq
{{i}\over{d-1}}\cdot\left(\log{{R}\over{\tilde
r}}+O_{\omega}({{1}\over{\tilde
r^{1/d-1}}}+O_\omega({{1}\over{R^\gamma}}))\right);$$ The involved
constants depend only on the point $p$, the metric $\omega$ and the
Evans kernel (and they are independent on the section $s$).
\endlemma

\Proof By \ref{propertiesofgp}, we can find a constant $a$ such
that, in $U_p$,  $\vert g_p\vert\Vert
z\Vert^{2(d-1)}\geq({{1}\over{d-1}}-a\Vert z\Vert^{\alpha})$. Thus,
The ball $B_z(t):=\{\Vert z\Vert^{2(d-1)}\leq {{1}\over
{t}}({{1}\over{d-1}}-{{a}\over{t^\alpha}})\}$ is contained in
$B_{g_p}(-t)$. As a consequence of [GK] Lemma 1.16 and Prop. 1.17 we
find that
$$\int_{div(s)_{g_p}(t)}\lambda^{d-1}\geq
\int_{div(s)\cap B_z(\vert t\vert)}\lambda^{d-1} \geq
{{i}\over{d-1}}\cdot{{1}\over{\vert t\vert}}\left(1-{{a}\over{\vert
t\vert^{\alpha}}}\right).$$

Because of our choices of coordinates on $U_p$,  in the open set
$g_p\leq -t$, we have that, $\omega^{d-1}\geq
\lambda^{d-1}\cdot\left(1+O_{\omega}({{1}\over{\vert
t\vert^{1/(d-1)}}})\right)$, where the involved constant depends
only on $\omega$. Consequently
$$\eqalign{\int_{-R}^{-\tilde r}dt\int_{div(s)_{g_p}(t)}\omega^{d-1}&\geq
\int_{-R}^{-\tilde r}dt\left(1+O_{\omega}({{1}\over{\vert
t\vert^{1/(d-1)}}})\right)\int_{div(s)_{g_p}(t)}\lambda^{d-1}\cr
&\geq \int_{-R}^{-\tilde r}{{i}\over{d-1}}\cdot {{1}\over{\vert
t\vert}}\cdot\left( 1+O_{\omega}({{1}\over{\vert
t\vert^{1/(d-1)}}})\right)\cdot\left(1-{{a}\over{\vert
t\vert^\alpha}}\right)dt.\cr}$$ The lemma follows.

We now prove the Theorem: By \ref{fmt} we have

$$\eqalign{T_{c_1(L)}(r)=&\lim_{R\to\infty}\int_{-R}^{\log(r)}dt\int_{B_{g_p}(t)}c_1(L)\wedge\omega^{d-1}=\cr
&\lim_{R\to\infty}\int_{-R}^{\log(r)}dt\int_{div(s)_{g_p}(t)}\omega^{d-1}
-\mu_p(\log\Vert s\Vert^2)(r)+\mu_p(\log\Vert
s\Vert^2)(e^{-R}).\cr}$$ By definition, of the norm on the strict
transform we have that
$$\mu_p(\log\Vert s\Vert)(e^{-R})=\int_{g_p=-R}\log\Vert \tilde
s\Vert^2d^cg_p\wedge\omega^{d-1}-{{i}\over{d-1}}\log(R).$$ We now
notice that
$$\int_{-R}^{\log(r)}dt\int_{div(s)_{g_p}(t)}\omega^{d-1}=\int_{-R}^{-\tilde r}dt\int_{div(s)_{g_p}(t)}\omega^{d-1}
+\int_{-\tilde r}^{\log(r)}dt\int_{div(s)_{g_p}(t)}\omega^{d-1};$$
consequently, if we let $R$ go to infinity and we apply the previous
lemma, the conclusion follows.
\endssection

\ssection Algebraic points on images of maps of finite order\par

Let $(A,\omega, g_p)$ be a conformally parabolic variety of
dimension $d$ with its Kh\"aler form and Evans kernel.

Suppose that $K$ is a number field (with a fixed embedding in $\CC$
as usual) and $X_K$ is a smooth projective variety of dimension $N$
defined over $K$. We fix an hermitian ample line bundle $L$ on $X_K$
and models etc. as in \S 2.

Let $S\subset X_K(\overline K)$ and $\alpha$ be a real number.

Let $\gamma: A\to X_K$ be an analytic map of finite order $\rho$.
For every point $P\in S\cap\gamma(A)$ we suppose that the germ of
$\gamma(A)$ is an $LG$ germ of type $\alpha$ defined over the
$K(P)$.

In this section we will prove that, under the condition above, we
can construct a closed positive $T$ current on $A$ with finite mass
and having Lelong number on each point of $\gamma^{-1}(S)$, bigger
then one. If one imagine this $T$ as the current associated to an
analytic divisor on $A$, the fact that its have finite mass is very
similar to an "algebraicity" condition while the condition on the
Lelong numbers corresponds to the fact that the divisor passes
trough $\gamma^{-1}(S)$. We will see that when $A$ is quasi
projective, this will imply that $\gamma^{-1}(S)$ is not Zariski
dense.

We will suppose that the hypotheses above are fixed once for all.

Let $S'\subset A$ be a countable set.

If $T$ is a current over $A$ and $P\in A$, we will denote by
$\nu(T,P)$ the Lelong number of $T$ in $P$.

\defi We will denote by $\Omega (S')$ the real number
$$\inf\left\{ \int_AT\wedge\omega^{d-1}\; /\; T\;{\rm is\; a\;
current\; of\; bidegree\; (1;1)\; with}\; \nu(T,P)\geq 1 \; \forall
P\in S'\right\}.$$ More generally, if $U$ is an open set of $A$, We
will denote by $\Omega(S';U)$ the real number
$$\inf\left\{
\int_UT\wedge\omega^{d-1}\; /\; T\;{\rm is\; a\; current\; on\; A\;
of\; bidegree\; (1;1)\; with}\; \nu(T,P)\geq 1 \; \forall P\in
S'\right\}.$$
\enddefi

A priori, $\Omega(S')$ may  be infinite. The aim of this section is
to show that, in the arithmetic situation, it is a finite number.

We give some tools to compute $\Omega(S')$.

\label Sfinite. proposition\par\prop Suppose that $S'=\bigcup_i S_i$
with $S_i\subseteq S_{i+1}$ and there exists a constant $M$ such
that $\Omega(S_i)\leq M$ for every $i$, then
$$\Omega(S')\leq M.$$
\endprop
\Proof For every positive $\epsilon$, we can find currents $T_i$
such that $\int_AT_i\wedge\omega^{d-1}\leq M$ and $\nu(T_i;P)\geq 1$
for every $P\in S_i$. By Ascoli--Arzela' Theorem we can find a
subsequence $T_{i_k}$ of the $T_i$ converging to a current $T$. By
construction $\int_AT\wedge\omega^{d-1}\leq M$ and, since
$\bigcup_{i_k}S_{i_k}=S'$, for every $P\in S$, we have $\nu(T;P)\geq
1$. The conclusion follows.
\smallskip
The proposition above is useful because it allows to suppose that
$S'$ has finite cardinality.

Let $\{ U_n\}_{n\in{{\Bbb N}}}$ be a sequence of relatively compact
open sets of $A$ such that

-- $\bigcup_nU_n=A$;

-- $U_n\subset U_{n+1}$.

We call such a sequence {\it an exhausting sequence}.

A proposition similar to \ref{Sfinite} allows to work with each term
of an exhausting sequence.

\label exhaustingS. proposition\par\prop Suppose that $\{
U_n\}_{n\in{{\Bbb N}}}$ is an exhausting sequence and that there
exists a constant $M$ such that, for every $n\in{\Bbb N}$ we have
$\Omega(S',U_n)\leq A$, then
$$\Omega(S')\leq M.$$
\endprop

\Proof For every positive $\epsilon$, and index $n$, we can find a
current $T_n$ such that $\int_{U_n}T_n\wedge\omega^{n-1}\leq
M+\epsilon$ and $\nu(T_n,P)\geq 1$ for every $P\in S$. By
Ascoli--Arzela' theorem, we can find a subsequence $\{ T_{1,1},
T_{1,2},\dots\}$ converging to a current $T^1$ on $U_1$ such that
$\int_{U_1}T^1\wedge\omega^{d-1}\leq M+\epsilon$ and $\nu(T;P)\geq
1$ for every $P\in S'\cap U_1$. We can extract from the subsequence
above a subsequence $\{ T_{2,1}, T_{2,2},\dots\}$ converging to a
current $T^2$ with the same properties on $U_2$ and so on. The
sequence $\{ T_{n,n}\}$ converges to a current $T$ with
$\int_{A}T^1\wedge\omega^{d-1}\leq M+\epsilon$ and $\nu(T;P)\geq 1$
for every $P\in S'$. The conclusion follows.
\smallskip
In the sequel, we fix an exhausting sequence $\{ U_n\}_{n\in{{\Bbb
N}}}$. Observe that such a sequence exist: it suffices to fix a
point $p\in A$ and take $U_n:=B_{g_p}(n)$.

Let $X$ be a projective variety and $L$ be an ample line bundle on
$X$ equipped with a positive metric (we also suppose that a metric
is fixed on $X$). Let $\gamma\colon A\to X$ be a map of finite order
$\rho$ with Zariski dense image. Let $S\subset A(\Bbb C)$.

Consider the linear map
$$\gamma_D:H^0(X,L^D)\hookrightarrow H^0(A,\gamma^\ast(L^D));$$
it is injective because the image of $\gamma$ is Zariski dense.

For every point $P\in S$, we denote by $I_P$ the ideal sheaf of $P$
in $A$. Given a finite subset $F\subset S$ we denote by $E_D^i$ the
kernel of the map obtained by composing $\gamma_D$ with the
restriction map $H^0(A, \gamma^\ast(L^D))\to\bigoplus_{P\in
F}(\O_A/I_P^i)\otimes L^D$ and by $\gamma_D^i$ the induced map
$$\gamma_D^i\colon E_D^i\longrightarrow\bigoplus_{P\in
F}S^i(T_PA^\ast)\otimes L^D_P.$$ Eventually, for every $P\in F$ we
denote by $\gamma_{D,P}^i$ the composite of $\gamma_D^i$ with the
canonical projection $\bigoplus_{P\in F}S^i(T_PA^\ast)\otimes
L^D_P\to S^i(T_PA^\ast)\otimes L^D_P$. We fix such a finite set $F$.

Since we fixed the  metric $\omega$ on $A$, all the vector spaces
involved are equipped with an hermitian metric.

\label firstboundonthenorms. theorem\par\thm Let $\epsilon>0$.
Suppose that we are in the situation above and $U\subset A$ is a
relatively compact open set. Then There exists a constant $C$
independent on $i$ and $D$ such that, for ${{i}\over{D}}\gg 0$ we
have
$$\log\Vert\gamma^i_{D,P}\Vert\leq
-i{{\Omega(F,U)}\over{\rho+\epsilon}}\log{{i}\over{D}}+C(i+D).$$
\endthm

\rmk One should compare this statement with Theorem
\ref{schwartzforparabolic}.
\endrmk

\Proof We fix $t_0\in\RR$ such that $B_{g_P}(t_0)\supset U$. Let
$s\in E_D^i$, then the current
$T_s:={{[div(\gamma^\ast(s))]}\over{i}}$ over $A$ is closed,
positive of bidegree $(1,1)$ and $\nu(T_s,Q)\geq 1$ for every $Q\in
F$. Thus
$$\int_UT_s\wedge\omega^{d-1}={{1}\over{i}}\int_{U\cap div(\gamma^\ast(s))}\omega^{d-1}\geq
\Omega(F,U).$$

Let $\tilde r$ as in \ref{generalboundofnonrms2}. For $\log(r)\geq
t_0$ we have
$$\eqalign{&\int_{-\tilde
r}^{\log(r)}dt\int_{div(\gamma^\ast(s))_{g_p}(t)}\omega^{d-1}\cr
&=\int_{-\tilde
r}^{t_0}dt\int_{div(\gamma^\ast(s))_{g_p}(t)}\omega^{d-1}+\int_{t_0}^{\log(r)}dt\int_{div(\gamma^\ast(s))_{g_p}(t)}\omega^{d-1}\cr
&\geq i\Omega(F,U)\left(\log(r)-t_0\right).\cr}$$ Consequently, by
theorem \ref{generalboundofnonrms2} and formula
\ref{boundonlinfinity} we find a constant, independent on $s$,  $i$
and $D$ such that
$$\log\Vert s\Vert^2_{\infty}-i\Omega(F,U)\log(r)+D\cdot T_\gamma(r)\geq
\log\Vert j^i(s)\Vert^2_{i,D}(P)+C(i+D).$$ Since $\gamma$ is of
finite order $\rho$, for every $\epsilon>0$ we can find $\lambda$
such that, for $r$ sufficiently big,  $T_\gamma(r)\leq\lambda
r^{\rho+\epsilon}$. The function $-i\Omega(F,U)\log(r)+D\lambda
r^{\rho+\epsilon}$ has a minimum when
$r^{\rho+\epsilon}=i{{\Omega(F;U)}\over{D\lambda(\rho+\epsilon)}}$;
thus there is a constant $C_1$ such that
$$\log\Vert
s\Vert^2_{\infty}-i{{\Omega(F;U)}\over{\rho+\epsilon}}\log{{i}\over{D}}\geq\log\Vert
j^i(s)\Vert^2_{i,D}(P)+C_1(i+D).$$ The conclusion follows.

\smallskip

Suppose that, as before, $K$ is a number field and
$\sigma:K\hookrightarrow\CC$ is an inclusion. We fix an algebraic
closure $\overline K$ of $K$ in $\CC$. As in the previous section,
We fix a smooth quasi projective variety $X_K$ over $K$ and
$S\subseteq X_K(\overline K)$. For every positive integer $r$ denote
by $S_r$ the set $\{ x\in S\; /\; [\QQ(x):K]\leq r\}$.

The main theorem of this section is the following.

\label generalbombieri. theorem\par\thm Suppose that $(A,\omega, g)$
is a conformally parabolic variety of dimension $d$. Let $\gamma:
A\to X_K(\CC)$ be an analytic map of finite order $\rho$ with
Zariski dense image. Suppose that, for every $P\in S\cap\gamma(A)$,
the formal germ $\hat A_P$ of $\gamma(A)$ near $P$ is (the pull
back) of an $LG$--germ of type $\alpha$ in its field of definition.
Then, for every $r$,
$${{\Omega(\gamma^{-1}(S_r))}\over{r}}\leq
{{N+1}\over{N-d}}\rho\alpha[K:\QQ].$$
\endthm
\Proof The proof is, {\it mutatis mutandis} identical to the proof
of theorem \ref{maintheoremforriemannsurfaces}. We give some
details. By propositions \ref{Sfinite} and \ref{exhaustingS} it
suffices to prove the following: Given a finite subset $F$ of
$\gamma^{-1}(S_r)$, a relatively compact set $U\subset A$  and
$\epsilon >0$ we have that
$${{\Omega(F,U)}\over{r}}\leq{{N+1}\over{N-d}}\alpha(\rho+\epsilon)[K:\QQ].$$
Let $F$ be such a set and $L$ be a finite extension of $K$ which
contains $\Bbb Q(\gamma(P_j))$ for every $P_j\in F$.

We  fix $P_j\in F$ and, and we fix an embedding $\sigma_{P_j}\in
S_{\Bbb Q(P_j)}$  extending $\sigma_0$. If $i/D\gg 0$, from
Corollary \ref{finitenorms}  and Theorem \ref{firstboundonthenorms}
we get,

$$\eqalign{{{1}\over{[\Bbb Q(P_j):\Bbb Q]}}&\sum_{\sigma\in M_{\Bbb
Q(P_j)}}\log\Vert \gamma_{D,j}^i\Vert_\sigma\cr &={{1}\over{[\Bbb
Q(P_j):\Bbb Q]}}\left( \log\Vert
\gamma_{D,j}^i\Vert_{\sigma_{\sigma_{P_j}}}+\sum_{\sigma\neq\sigma_{P_j}}\log\Vert
\gamma_{D,j}^i\Vert_\sigma\right)\cr &\leq
-{{\Omega(F;U)}\over{(\rho+\epsilon)\cdot[\Bbb Q(P_j):\Bbb
Q]}}i\log\left({{i}\over{D}}\right)+\alpha\cdot i\log
(i)+C_2\left(i+D)\right),\cr}$$ with $C_2$ independent on $i$, $D$
and $\Bbb Q(P_j)$. Thus, we deduce
$$\eqalign{{{1}\over{[L:\Bbb Q]}}&\sum_{\sigma\in
M_L}\log\Vert\gamma_D^i\Vert_\sigma\leq{{1}\over{[L:\Bbb
Q]}}\sum_{\sigma\in M_L}\sup_{P_j\in F}\log\Vert
\gamma_{D,j}^i\Vert_\sigma+\log(Card(F))\cr &\leq
-{{\Omega(F;U)}\over{r\cdot[K:\Bbb
Q](\rho+\epsilon)}}i\log\left({{i}\over{D}}\right)+\alpha\cdot i\log
(i)+C_2\left(i+D\right)+\log(Card(F)).\cr}$$ But since the image is
Zariski dense, we can apply Theorem \ref{generalalgebraization} and
we obtain that
$${{\Omega(F;U)}\over{r\cdot[K:\Bbb Q]}}\leq {{N+1}\over{N-d}}(\rho+\epsilon)\alpha.$$
The conclusion follows.

\smallskip

By taking weak limits of currents one obtain:

\label generalbombieri2. theorem\par\thm Suppose that we are in the
hypotheses as above. Then there exists a closed positive current $T$
of bidegree $(1,1)$ over $A$ such that:

-- $\int_A T\wedge \omega^{d-1}\leq
2{{N+1}\over{N-d}}\rho\alpha[K:\QQ]$;

-- for every $P\in\gamma^{-1}(S)$ we have that
$\nu(T;P)\geq{{1}\over{[K(\gamma(P)):K]}}$.
\endthm

\Proof By  theorem \ref{generalbombieri}, for every integer $r$ and
$\epsilon>0$ we can find a current $T_r$ such that
${{1}\over{r}}\int_A T_r\wedge \omega^{d-1}\leq
{{N+1}\over{N-d}}\rho\alpha[K:\QQ]+\epsilon$ and $\nu(T_r; P)\geq 1$
for every $P\in\gamma^{-1}(S_r)$.

For every positive integer $R$ consider the current
$T_R:=\sum_{r=1}^R{{T_r}\over{r^2}}$.

By construction $\int_AT_R\wedge\omega^{d-1}\leq \left(
\sum_{r=0}^R{{1}\over{r^2}}\right)\left({{N+1}\over{N-d}}\rho\alpha[K:\QQ]+\epsilon\right)$.
Thus we may find a subsequence of the $\{ T_R\}$ converging to a
current $T$ such that $\int_AT\wedge\omega^{d-1}\leq \left(
\sum_{r=0}^\infty{{1}\over{r^2}}\right)\left({{N+1}\over{N-d}}\rho\alpha[K:\QQ]+\epsilon\right)\leq
2\left({{N+1}\over{N-d}}\rho\alpha[K:\QQ]+\epsilon\right)$.

Let $P\in \gamma^{-1}(S_r)\setminus\gamma^{-1}(S_{r-1})$. For every
$i\geq r$, the current ${{T_{i}}\over{i^2}}$ has Lelong number on
$P$ bigger or equal then ${{1}\over{i^2}}$. Consequently the Lelong
number of $T$ on $P$ is bigger or equal then
$\sum_{i=r}^{\infty}{{1}\over{i^2}}$. Since
$\sum_{i=r}^{\infty}{{1}\over{i^2}}\geq\int_r^\infty{{1}\over{t^2}}dt={{1}\over{r}}$.
the conclusion follows

\rmk One can obtain as corollary a statement analogous to
\ref{bostgk} (but one needs to change a little bit the proof).
However it will be again a consequence of [Bo], thus we leave the
details to the reader.\endrmk

\endssection

\ssection Rational points on affine varieties\par In this subsection
we will show some consequence of our theory for affine varieties.

Often one has an analytic map between an affine and a projective
variety and one want to show that, if the map is not algebraic, then
the preimage of the algebraic points is not Zariski dense under some
condition. We will show that theorem \ref{generalbombieri} will
easily imply this for maps of finite order (and $LG$ germs). One
should notice that theorem \ref{generalbombieri} is more general and
may have other applications.

Before we prove the theorem, we want to show that the value
distribution theory we developed here is essentially equivalent to
the value distribution theory developed by Griffiths and King in
[GK].

We begin by fixing the hypotheses and the notations as in [GK].

Let $A$ be a smooth affine variety of dimension $d$ defined over
$\Bbb C$. We can embed $A$ as a Zariski open set
$A\hookrightarrow\overline A$, where $\overline A$ is a smooth
projective variety. Moreover we can suppose that the closed set at
infinity, $D_{\infty}:=\overline A\setminus A$,  is a divisor whose
support is a divisor with simple normal crossing.

We fix such a compactification of $A$ and a projective embedding
$\iota\colon\overline A\hookrightarrow\Bbb P^M$ for a suitable $M$.
We can suppose that there is an hyperplane $H_{\infty}$ in $\Bbb
P^M$ such that $D_{\infty}=H_{\infty}\cap\overline A$ (without
multiplicity). We will denote by $R$ the degree of $\overline A$ in
$\Bbb P^M$. We will denote by $\O_{\UA}(1)$ the restriction of the
ample line bundle $\O_{\Bbb P^M}(1)$ to $\UA$; we will suppose that
it is equipped with the Fubini--Study metric $\Vert\cdot\Vert_{FS}$.
To simplify notations we will suppose that $H^i(\UA;\O_{\UA}(1))=0$
for every $i>0$ (this hypothesis is not indispensable, nevertheless
it can always be supposed and simplify the statement of the main
result).

Let $\omega$ be the standard K\"ahler form on $\overline A$ induced
from the embedding $\iota$.

Let $p\in A$ be a point.

We will now recall the classical definition of counting function,
given for instance in [GK].

Choose a linear subspace $L\simeq\Bbb P^{M-(d+1)}\subset H_{\infty}$
which do not meet $\overline A$ and we project from $L$ on a
suitable projective space $\Bbb P(V)\simeq\Bbb P^d$.

We get a commutative diagram
$$\matrix{\overline {A}&\buildrel{\iota}\over\llongrightarrow&{\Bbb P}^M\setminus L\llongleftarrow &H_{\infty}\setminus L\cr
    \mapdownl{\pi}& &\mapdownl{\pi} &\Big\downarrow\cr
    {\Bbb P^d} & = &{\Bbb P^d} \llongleftarrow &{\Bbb P}^{d-1}\cr}$$
with $\pi^{-1}(\Bbb P^{d-1})\cap\overline A=H_{\infty}$. From this
we obtain a finite branched covering $\pi\colon A\to\Bbb C^d$.
Denote by $Ram(\pi)$ the ramification divisor of $\pi$. We suppose
the situation above fixed once for all.

We fix some coordinate on $\Bbb C^d$ and we define $\varphi
:=\pi^{\ast}(\log\Vert z\Vert^2)$, where, for $z=(z_1,\dots,
z_d)\in\Bbb C^d$, we define $\Vert z\Vert^2:=\sum_i\vert z\vert^2$.

If $\alpha$ is a $(1,1)$ form on $A$ we define the Nevanlinna
counting function with respect to $\omega$ to be
$$\eqalign{\tilde T_{\alpha}(r)&:=\int_0^r{{dt}\over{t}}\int_{B_\varphi(\log(t))}\alpha\wedge(dd^c\varphi)^{d-1}\cr
&=\int_{-\infty}^{\log(r)}dt\int_{\varphi\leq
t}\alpha\wedge(dd^c\varphi)^{d-1}\cr &
=\int_A(\log(r)-\varphi)^+(\alpha\wedge dd^c\varphi)^{d-1}\cr}$$
(where, as in \S 4, $(f)^+$ means $\sup\{ f, 0\}$).

Let $X$ be projective variety and $\gamma\colon A\to X$ be an
analytic map. Suppose that $\overline L$ is an ample line bundle
equipped with a positive metric. We define $\tilde T_{\gamma}(r)$ to
be $\tilde T_{\gamma^{\ast}(c_1(\overline L))}(r)$. In [GK] chapter
2 many properties of $\tilde T_{\gamma}(r)$ are proved (in
particular the fact that it is, up to bonded functions, independent
on the chosen metric).

\label classfiniteorderofaffine. definition\par\defi We say that
$\gamma\colon A\to X$ is {\it of finite Nevanlinna order}
$\tilde\rho$ if
$$\limsup_{r\to\infty}{{\log(\tilde T_\gamma(r))}\over{\log(r)}}=\tilde\rho$$
\enddefi

One sees that, when $A$ is an affine variety, we have two possible
value distribution theory on it.  We will show that these two
theories are comparable, namely the two definitions of order of an
analytic map coincide:

\label comparisonoforders. proposition\par\prop If $\gamma:A\to X$
is an analytic map then
$$\rho=\tilde\rho.$$
\endprop

We first need three lemmas. Although they are standard and easy, we
provide a proof for reader's convenience.

If $L\colon \Bbb R\to \Bbb R$ is a positive increasing function, we
will denote by $\rho_L$ the number
$\limsup_{r\to\infty}{{\log(L(r))}\over{\log(r)}}.$

\label lang1. lemma\par\lemma Let $E$ be a subset of $\Bbb R$ having
finite Lebesgue measure. Then, for $R$ sufficiently big,
$$\rho_L=\limsup_{r\geq R,r\not\in E}{{\log(L(r))}\over{\log(r)}}:=\rho_E.$$
\endlemma
\Proof It is evident that $\rho_L\geq \rho_E$. Suppose that
$\rho_L\geq\rho_E+\epsilon$. Let $S\gg 0$ such that $meas(E\cap
[S;\infty])\leq\delta$ (we will fix delta at the end of the proof).
There are infinitely many $r\geq S$ such that $\log(L(r))>(\rho_E
+\epsilon)\log(r)$. Fix one of them. Then $(r,
r+\delta)\not\subseteq E$, thus there exists $s\in (r,
r+\delta)\setminus E$. Consequently
$$(\rho_E+\epsilon)\log(r)\leq \log(L(r))\leq \log(L(s))\leq
(\rho_E+\epsilon/2)\log(s)\leq
(\rho_E+\epsilon/2)(\log(r)+O({{\delta}\over{r}}))$$ and this is a
contradiction.

\label lang2. lemma\par\lemma (Lang) Let $\psi$ be a positive
function such that
$\int_{1+\epsilon}^{\infty}{{1}\over{x\psi(x)}}dx<\infty$. Let $T$
be positive increasing function. Then there exists a subset
$E\subset\Bbb R$ with
$meas(E)\leq\int_{1+\epsilon}^{\infty}{{1}\over{x\psi(x)}}dx<\infty$,
such that, for every $x\not\in E$,
$$T'(x)\leq T(x)\psi(T(x)).$$
\endlemma

\Proof Consider the subset $F:=\{ x \;\; {\rm s.t}\;\; T'(x)\geq
T(x)\psi(T(x))\}$. If we denote by ${\Bbb I}_F(x)$ the
characteristic function of $F$, then one easily see that $${\Bbb
I}_F(x)\leq {{T'(x)}\over{T(x)\psi(T(x))}}.$$ Integrating the
inequality above and making the change of variables $T(x)=y$ the
conclusion follows.

\label usefulang2. corollary\par\cor Suppose that $T$ is an
increasing derivable function, then $$\rho_{T'}\leq\rho_T.$$
\endcor
\Proof Take  $\psi(x)=\log^{1+\epsilon}(x)$ and apply the two
lemmas.

\label comparitionofforms. lemma\par\lemma Let $\omega$ be the
K\"ahler form on $\UA$ then, on $\UA\setminus\pi^{-1}(0)$,
$$2\omega\geq dd^c(\varphi).$$
\endlemma
\Proof Since everything is invariant under the action of the unitary
group $U(d+1)\times U(M-d)$, it suffices to verify the inequality on
the point with homogeneus coordinates $[1,0,\dots,0,1,0,\dots, 0]$
(the first sequence is a sequence of $d$ zeros). On that point an
explicit calculation gives the result.

We can now start the proof of \ref{comparisonoforders}.

\Proof (of \ref{comparisonoforders}). By \ref{comparitionofforms}
one easily see that $\rho\geq \tilde\rho$. We need to prove the
converse inequality.

Take $a\in\Bbb R$ sufficiently small, denote by $T_{\gamma}^1(r)$
(resp. $\tilde T_{\gamma}^1(r)$) the function
$$\int_{a}^{\log(r)}dt\int_{a\leq g_p\leq
t}c_1(L)\wedge\omega^{d-1}$$ (resp.
$\int_{a}^{\log(r)}dt\int_{a\leq\varphi\leq t}c_1(L)\wedge
(dd^c(\varphi)^{d-1}$) then
$\rho=\limsup{{\log(T_{\gamma}^1(r))}\over{\log(r)}}$ (resp.
$\tilde\rho=\dots$).

Moreover, denoting by $T^2_{\gamma}(r)$ the function
$$\int_{a}^{\log(r)}dt\int_{a\leq \varphi\leq
t}c_1(L)\wedge\omega^{n-1},$$ we have
$\rho=\limsup{{\log(T_{\gamma}^2(r))}\over{\log(r)}}$. Indeed, there
exists a $C^{\infty}$ function $f$ on $A\setminus B_{g_p}(a)$ such
that $g_p={{1}\over{\deg(\overline A)}}\varphi+f$; indeed, on
$A\setminus B_{g_p}(a)$, we have that
$dd^c(\varphi)=-\delta_{H_{\infty}}+\alpha$, where $\alpha$ is a
closed positive (1,1) form; thus
$dd^c(g_p-\varphi)\wedge\omega^{d-1}$ is a $C^{\infty}$ form.

Denote by $\theta$ the positive (1,1) form $dd^c(\varphi)$ and
$$T^p(r):=\int_{a}^{\log(r)}dt\int_{a\leq \varphi\leq
t}c_1(L)\wedge\theta^{(n-1)-p}\wedge\omega^{p}.$$ By induction, it
suffices to prove that $\rho_{T^p}=\rho_{T^{p-1}}$.

It is evident that $\rho_{T^p}\geq\rho_{T^{p-1}}$.

Denote the coordinates on $\Bbb C^d$ by $(z_1,\dots, z_d)$ and by
$u(z)$ the function $\pi^{\ast}(\sum\vert z_i\vert^2)$.  A simple,
computation, obtained by taking successive projections on
hyperplanes of codimension one, shows that $\omega-\theta=dd^c(k)$
where $k$ is the function $\log(1+{{1}\over{h}}+k_1)$, where $k_1$
is a function which is bounded and extends to a $C^{\infty}$
function on $\overline A\setminus B_{g_p}(a)$. Observe that,
consequently,  $k$ is bounded on $A\setminus B_{g_p}(a)$.

We now compute, using Stokes theorem,
$$\eqalign{&T^{p+1}(r)-T^{p}(r)=\int_{a}^{\log(r)}dt\int_{a\leq \varphi\leq
t}c_1(L)\wedge(\theta^{(d-1)-(p+1)}\wedge\omega^{p+1})-\theta^{(d-1)-p}\wedge\omega^{p}\cr
&=\int_{a}^{\log(r)}dt\int_{a\leq \varphi\leq
t}c_1(L)\wedge\theta^{(d-1)-(p+1)}\wedge\omega^{p}\wedge(dd^c(k))\cr
&=\int_{\varphi=\log(r)}kd^c(\varphi)\wedge\theta^{(d-1)-(p+1)}\wedge\omega^{p}-\int_{a\leq\varphi\leq\log(r)}k
dd^c(\varphi)\wedge\theta^{(d-1)-(p+1)}\wedge\omega^{p}+ C.\cr}$$
But, since $k$ is bounded, we apply Stokes theorem again and find a
constant $C$ such that
$$T^{p+1}(r)-T^{p}(r)\leq C\cdot
\int_{a\leq\varphi\leq\log(r)}c_1(F)\wedge\theta^{(d-1)-p}\wedge\omega^p={{dT^{p}}\over{dr}}(r).$$
The conclusion follows from \ref{usefulang2}.

\label omegam. definition\par\defi Let $S\subset A(\Bbb C)$; we will
say that $\varpi(S)=m$ if the natural map
$$\beta_m:H^0(\UA;\O_{\UA}(m)\longrightarrow \prod_{p\in S}\O_{\overline
A}(m)\vert_p$$ {\rm is  not} injective but the map
$$\beta_{m-1}\colon H^0(\UA;\O_{\UA}(m-1)\longrightarrow \prod_{p\in S}\O_{\overline
A}(m-1)\vert_p$$ is injective.
\enddefi
\rmk a) If $\beta_m$ is not injective, then, for every $\ell\geq 0$,
$\beta_{m+\ell}$ is not injective.

b) If $\varpi(S)=m$ then $S$ is contained in the pull back of a
divisor of degree $m$ on $\Bbb P^M$ and {\it is not} contained in
any pull back of divisors of degree strictly less then $m$.
\endrmk

Let $X_K$ be a projective variety of dimension $N$ and $S\subseteq
X_K(\overline K)$; We will define $S_r$ as before.

\label affinebombieri. theorem\par\thm Let $\gamma:A\to X_K(\CC)$ be
an analytic map of finite order $\rho$. Suppose that, for every
$P\in S\cap\gamma(A)$, the germ of $\gamma(A)$ near $P$ is
isomorphic to  an $LG$ germ of type $\alpha$ in the field of
definition of $P$. Then for every positive integer $r$ we have
$${{\varpi(\gamma^{-1}(S_r))}\over{r}}\leq{{N+1}\over{N-d}}{{d!\rho\cdot\alpha}\over{\deg(\UA)}}.$$
\endthm

As a consequence, for every $r$, $\gamma^{-1}(S_r)$ is {\it not}
Zariski dense in $A$.

\Proof Denote by $V_m$ the vector space $H^0(\UA,\O_{\UA}(m))$ and
fix a basis $\{ s_0;s_1;\dots; s_\ell\}$ of $V_1$ such that
$div(s_0)=D_\infty$. Suppose that $\varpi(\gamma^{-1}(S_r))>m$. We
can find a finite set $F\subset\gamma^{-1}(S_r)$ such that the map
$$\beta_m\colon V_m\longrightarrow \prod_{P\in F}\O_{\overline
A}(m)\vert_P$$ is an isomorphism. Moreover, since $s_0(P)\neq 0$ for
every $P\in F$, we have a commutative diagram
$$\matrix{V_m&\llongmaprighto{\simeq}&\prod_{P\in F}\O_{\overline
A}(m)\vert_P\cr \mapdownr{\cdot s_0} & & {\simeq} \mapdownr{\cdot
s_0}\cr V_{m+1}&\llongrightarrow & \prod_{P\in F}\O_{\overline
A}(m+1)\vert_P.\cr}$$

For every multindex $I=(i_0,i_1,\dots, i_\ell)$ we denote by $\vert
I\vert$ the sum $i_0+i_1+\dots+i_\ell$ and by $s^I$ the section
$s_0^{i_0}s_1^{i_1}\cdots s_\ell^{i_\ell}\in V_{\vert I\vert}$.
Because of the commutative diagram above, for every multindex $I$
with $\vert I\vert=m+1$ and $i_0=0$ we can find an element $R_I\in
V_m$ such that $\beta_{m+1}(s_I+R_I\cdot s_0)=0$. We will denote by
$Q_I$ the section $s_I+R_I\cdot s_0$.

The linear system generated by the $Q_I$ and $s_0^{m+1}$ define a
map $Q\colon \UA\to \Bbb P^h$.

We claim that the map $Q$ is defined everywhere and finite: indeed,
it is defined everywhere because, $s_0$ do not vanish outside the
hyperplane at infinity $D_\infty$ and, the $Q_I$ cut over this
hyperplane the complete linear system $H^0(\O(m+1))$. It is finite:
we first observe that the restriction of $Q$ to the hyperplane at
infinity is an embedding (it is the restriction of the $m+1$--th
Veronese embedding). There is an hyperplane $H$ in $\Bbb P^h$, the
pull back of which is the hyperplane at infinity $D_\infty$; if
there was a fibre of positive dimension this would cut the
hyperplane at infinity somewhere and this is not possible.

We have that $Q(F)=[0:0:\dots:0:1]$.

Let $b:\tilde A\to\UA$ be the blow up of $\UA$ in
$Q^\ast([0:0:\dots:0:1])$ and $E$ be the exceptional divisor. By
projection, we have a map
$$\tilde P:\tilde A\longrightarrow \PP^{h-1}$$
and $\tilde p^{\ast}(O(1))=b^\ast(\O_{\UA}(m))-E$. Thus
$b^\ast(\O_{\UA}(m))-E$ is a nef line bundle on $\tilde A$.

By theorem \ref{generalbombieri} we can find a closed positive
current $T$ on $A$ of bidegree $(1,1)$ such that

-- ${{1}\over{r}}\int_AT\wedge \omega^{d-1}\leq
{{N+1}\over{N-d}}\rho\alpha[K:\QQ]$;

-- For every $P\in F$, we have $\nu(T,P)\geq 1$.

Let $U:=A\setminus \{ F\}$ and $\II_U$ its characteristic function.
By Skoda--El Mir Theorem, the current $\II_UT$ extends to a closed
positive current $\tilde T$ on $\tilde A$. By [GK] lemma 1.16 and
prop. 1.17 (to be precise by a small variation of them, cf. remark
after lemma 1.16 of loc cit.), and the fact that $\nu(T,P)\geq 1$
for every $P\in F$ we obtain
$$(\tilde T; E^{d-1})\geq Card(F).$$
Since $\tilde T$ is closed and positive and $b^\ast(\O_{\UA}(m))-E$
is nef,
$$(\tilde T; (b^\ast(\O_{\UA}(m))-E)^{d-1})=m^{d-1}(\tilde
T;(b^\ast(\O_{\UA}(1)^{d-1})))-(\tilde T; E^{d-1})\geq 0.$$ Thus
$$m^{d-1}r{{N+1}\over{N-d}}\rho\alpha[K:\QQ]\geq
Card(F)=h^0(\UA;\O_{\UA}(m))={{m^d}\over{d!}}\deg(\UA).$$ The
conclusion follows.

\smallskip

We remark the following corollary where the fact that Theorem
\ref{affinebombieri} is a generalization of the Bombieri Schneider
Lang criterion become evident:

\label classicalbsl. corollary\par\cor  Let $X$ be a projective
variety defined over a number field $K$ and let $F\hookrightarrow
T_X$ be a foliation. Suppose that $\gamma\colon A\to X$ is an
analytic map of finite order $\rho$ from an affine variety to $X$
such that the image is a non algebraic leaf of the foliation. Then
$\gamma^{-1}(X(K)\setminus Sing(F))$ is not Zariski dense.\endcor

We find Bombieri Schneider Lang Theorem by taking $X=\Bbb C^N$ and
$A=\Bbb C^n$.

\endsection

\

\section References\par

\

\item{[AS]} L.V. Ahlfors and L. Sario, {\it Riemann Surfaces},
Mathematical Series, {\bf 26}, Princeton, N.J., Princeton University
Press, (1960).

\item{[Au]} T. Aubin. {\it Nonlinear analysis on manifolds. Monge-Ampère equations}.
Grundlehren der Mathematischen Wissenschaften, 252. Springer-Verlag,
New York, 1982

\item{[Bom]} E. Bombieri, {\it Algebraic values of meromorphic
maps,} Inv. Math. {\bf 10}, (1970) pp. 267--287.

\item{[Bo]} J.B. Bost, {\it Algebraic leaves of algebraic
foliations over number fields}, Publ. Math, I.H.E.S. {\bf 93},
(2001) pp.161--221.

\item{[Bo2]} J.B. Bost: Private mail to the author.

\item{[Dl]} P. Deligne, {\it Equations differentielles \'a points
singuliers r\'eguliers}, Springer Lect. N. in Math. 163, Springer
Verlag, 1970.

\item{[De1]} J.P. Demailly, {\it Formule de Jensen en plusieur
variables et applications arithmetiques}, Bull. Soc. Math. France,
{\bf 110}, (1982) pp. 75--102.

\item{[De2]} J.P. Demailly, {\it Mesures de Monge--Ampere et
mesures pluriharmoniques}, Math. Zeit., {\bf 194}, (1987) pp.
519--564.

\item{[De3]} J.P. Demailly, {\it Complex analytic and differential
geometry}, Book preprint available at
http://www-fourier.ujf-grenoble.fr/~demailly/books.html.

\item{[GK]} P.A. Griffiths and J. King, {\it Nevanlinna theory and
holomorphic mappings between algebraic varieties}, Acta Math., {\bf
130}, (1973) pp. 145--220.

\item{[Gr]} A. Grygor'yan, {\it Analytic and geometric background
of recurence and non--ex\-plo\-sion of the Brownian motion on
Riemaniann manifolds}, Bull. Amer. Math. Soc. (N.S.), {\bf 36},
(1999), pp. 135--249.

\item{[La]} S. Lang, {\it Introduction to transcendental numbers},
Addison--Wesley Publishing Co. Reading, Mass.--London--Don Mills,
Ont. (1966).

\item{[MQ]} M. McQuillan, {\it Non commutative Mori theory}, Preprint I.H.E.S. , 141 pp.

\item{[SN]} L. Sario and M. Nakai, {\it Classification theory of
Riemann surfaces}. Die Grun\-dle\-hren der ma\-the\-ma\-ti\-schen
Wissenschaften, Band 164 Springer-Verlag, New York-Berlin 1970.

\item{[Vo]}  C. Voisin. {\it Hodge theory and complex algebraic geometry}.
Cambridge Studies in Advanced Mathematics, 76. Cambridge University
Press, Cambridge, 2002.

\item{[Wa]} F. Warner. {\it Foundations of differentiable manifolds and Lie
groups.} Graduate Texts in Mathematics, 94. Springer-Verlag, New
York-Berlin, 1983.

\endsection

\bigskip

C.~Gasbarri:  Dipartimento di Matematica dell'Universit\`a di Roma
``Tor Vergata", Viale della Ricerca Scientifica, 00133 Roma (I).

gasbarri@mat.uniroma2.it

\end